\documentclass[11pt]{amsart}
\usepackage{amsmath}
\usepackage{amsxtra}
\usepackage{amscd}
\usepackage{amsthm}
\usepackage{amsfonts}
\usepackage{amssymb}
\usepackage{eucal}
\newcommand{\g}{{\mathfrak{g}}}
\newcommand{\gln}{{\mathfrak{gl}}_n}
\newcommand{\Sym}{{\mathfrak{S}}}
\newcommand{\slt}{\mathfrak{sl}_2}
\newcommand{\slth}{\widehat{\mathfrak{sl}}_2}

\newcommand{\sln}{\mathfrak{sl}_n}
\newcommand{\nn}{\nonumber}
\newcommand{\bea}{\begin{eqnarray}}
\newcommand{\ena}{\end{eqnarray}}
\newcommand{\be}{\begin{eqnarray*}}
\newcommand{\en}{\end{eqnarray*}}

\newcommand{\ch}{{\rm ch}}
\newcommand{\End}{\mathop{{\rm End}}}

\newcommand{\gr}{\mathop{{\rm gr}}}
\newcommand{\br}[1]{{\langle #1 \rangle}}  
\newcommand{\qbin}[2]
{{ 
\left[
\begin{matrix}{\displaystyle #1}\\
{\displaystyle #2}\end{matrix}
\right]
}}

\numberwithin{equation}{section}
\newtheorem{thm}{Theorem}[section]
\newtheorem{prop}[thm]{Proposition}
\newtheorem{lem}[thm]{Lemma}
\newtheorem{cor}[thm]{Corollary}

\newtheorem{defn}{Definition}

\renewcommand{\cong}{\simeq}

\newcommand{\gc}{{\g[t]}}
\newcommand{\zz}{{\mathcal Z}}
\newcommand{\fusn}{\circledast}
\newcommand{\wg}{\widehat{\g}}

\newcommand{\simk}{{P_+^{(k)}}}
\newcommand{\ant}{S}
\newcommand{\hkr}{{\theta^\vee}}
\newcommand{\fil}{{\mathcal F}_\zz}
\newcommand{\gn}{{\mathfrak n}}
\newcommand{\hh}{{\mathfrak h}}
\newcommand{\ann}{{\rm Ann}}
\newcommand{\gX}{{\mathfrak X}}
\newcommand{\gY}{{\mathfrak Y}}
\newcommand{\ik}{{\mathfrak I}^{(k)}}
\newcommand{\tx}{\tilde{X}}
\newcommand{\txx}{\tilde{\gX}}
\newcommand{\fk}{{\mathcal K}}
\newcommand{\ga}{{\mathfrak a}}
\newcommand{\cpu}{U_{\{0\}}(\gc)}
\newcommand{\Ikla}{I^{(k)}_{\lambda}}
\newcommand{\Ila}{I_{\lambda}}
\newcommand{\Lkla}{L^{(k)}_{\lambda}}
\newcommand{\vkla}{v^{(k)}_{\lambda}}
\newcommand{\Lkinf}{L^{(k)}_{\lambda}(\infty)}
\newcommand{\vkinf}{v^{(k)}_{\lambda}(\infty)}
\newcommand{\C}{{\mathbb C}}
\newcommand{\Z}{{\mathbb Z}}

\newcommand{\F}{\mathcal F}
\newcommand{\Cc}{\mathcal C}
\newcommand{\Lc}{\mathcal L}
\newcommand{\tg}{\tilde \g}

\begin{document}
\title[Coinvariants and fusion product]
{Spaces of coinvariants and fusion product I.\\
{}From equivalence theorem to Kostka polynomials}
\author{B. Feigin, M. Jimbo, R. Kedem, S. Loktev, T. Miwa}
\address{BF: Landau institute for Theoretical Physics, Chernogolovka,
142432, Russia}\email{feigin@feigin.mccme.ru}  
\address{MJ: Graduate School of Mathematical Sciences, The 
University of Tokyo, Tokyo 153-8914, Japan}\email{jimbomic@ms.u-tokyo.ac.jp}
\address{RK: Department of Mathematics, University of Illinois,
1409 W. Green St. Urbana, IL 61801, USA}
\email{rinat@math.uiuc.edu}
\address{SL:  Institute for Theoretical and Experimental Physics,
  B. Cheremushkinskaja, 25, Moscow 117259, Russia}
\address{
 Independent University of Moscow, B. Vlasievsky per, 11, 
Moscow 121002, Russia}\email{loktev@mccme.ru}
\address{TM: Division of Mathematics, Graduate School of Science, 
Kyoto University, Kyoto 606-8502
Japan}\email{tetsuji@kusm.kyoto-u.ac.jp}
\date{\today}
\begin{abstract}
The fusion rule gives the dimensions of spaces of conformal blocks
in the WZW theory. We prove a dimension formula similar to the fusion rule
for spaces of coinvariants of affine Lie algebras $\widehat\g$.
An equivalence of filtered spaces is established between spaces
of coinvariants of two objects: 
highest weight $\widehat\g$-modules and 
tensor products 
of finite-dimensional evaluation representations of $\g\otimes\C[t]$.

In the $\slth$ case we prove that their associated graded spaces
are isomorphic to the spaces of coinvariants of fusion products, 
and that 
their Hilbert polynomials are the level-restricted Kostka polynomials.
\end{abstract} 
\maketitle 

\setcounter{section}{0}
\setcounter{equation}{0}

\section{Introduction}\label{sec:1}
First recall the notion of conformal blocks in the Wess-Zumino-Witten
conformal field theory. Let $\g$ be a simple Lie algebra over $\C$, let
$\tg=\g\otimes\C[t,t^{-1}]$ be the loop algebra, and $\tg^{\{0\}}$ its
completion with respect to power
series in $t$. As a geometric part of the data consider a compact complex 
curve $C$ with chosen distinct points $p_0, \dots, p_N$. We take these points
to be regular, and choose local coordinates $t_0, \dots, t_N$. 
By $\g^{C\setminus \{p_0, \dots, p_N\}}$
denote the Lie algebra of $\g$-valued meromorphic functions on $C$
with possible poles at $p_i$. Local coordinates provide us with 
the inclusion
\begin{equation}\label{LIE}
\g^{C\setminus\{p_0,\dots,p_N\}}\to\tg^{\{0\}}\oplus\dots\oplus\tg^{\{0\}}
\oplus\C K
\end{equation}
sending a function into its Laurent expansion at each point. Here $\C K$
is the one-dimensional center. The central extension in the right hand side
is given by the sum of canonical cocycles of each summand $\tg^{\{0\}}$.
It is trivial on the image of $\g^{C\setminus\{p_0,\dots,p_N\}}$.

Let $L_0,\dots, L_N$ be irreducible highest weight representations
of a fixed level $k$ of $\wg^{\{0\}}=\tg^{\{0\}}\oplus\C K\oplus\C d$.
Let $v_i$ be their highest weight vectors.
\footnote{In the main text, we will consider reducible highest weight
representations with cyclic vectors $v_i$. In this introduction, we restrict to
irreducible representations in order to ease the presentation.}

Define the space of conformal blocks by 
\begin{equation}\label{CB}
\left<L_0,\dots,L_N\right>^C_{p_0,\dots,p_N}=
L_0\boxtimes\dots\boxtimes L_N/ \g^{C\setminus \{p_0, \dots, p_N\}}
(L_0\boxtimes\dots\boxtimes L_N).
\end{equation}
In this paper, we call it the space of conformal coinvariants.
Here we use the symbol $\boxtimes$ in place of the tensor product symbol
$\otimes$ in order to emphasize that we consider the action of the
Lie algebra (\ref{LIE}).

If $L_i$ are integrable, then it was proved by \cite{TUY}
that these spaces are finite dimensional, 
and that their dimensions are given by 
the fusion rule, in particular, it depends only on representations
and the genus of the curve (see Theorem \ref{thm:verlinde}
for the fusion rule when the curve is $\C P^1$).

Introduce a filtration on the space of conformal coinvariants in the
following way. Let $v_1, \dots, v_N$ be the highest weight vectors of
representations $L_1, \dots, L_N$, respectively. Then we have a map from
$L_0$ to the space of conformal coinvariants given by
\bea\label{map}
L_0\simeq L_0\otimes v_1\otimes\dots\otimes v_N\to
\left<L_0,\dots,L_N\right>^C_{p_0,\dots,p_N}.
\ena
By a standard argument (see Subsection \ref{subsec:2.4}) 
one can show that the map~(\ref{map}) is surjective.
This gives us a realization of the space of conformal coinvariants as a
quotient of $L_0$. Since the representation  $L_0$ is graded by degree of
$t_0^{-1}$, our space inherits the corresponding increasing filtration. 

In this paper, we study this filtration for rational $C$ and
irreducible highest weight representations $L_0,\dots,L_N$
of level $k\in\Z_{>0}$.

{}From this point we restrict to the case of rational $C$.
We take $p_0=\infty$, and $p_i=z_i\in\C$ $(1\leq i\leq N)$.
We use the local coordinates $t_0=1/t$ and $t_i=t-z_i$. We take $L_0$ to be
the level $k$ irreducible highest weight module whose degree zero subspace
is the irreducible $\g$-module $\pi_\lambda$ with highest weight $\lambda$.
We denote this $L_0$ specifically by $L^{(k)}_\lambda(\infty)$.
Similarly, we take highest weights $\lambda_i$ for $L_i$ $(1\leq i\leq N)$.

The space $\Lkinf$ is generated from the highest weight
vector $v_0$ by the action of $\g[t]=\g\otimes\C[t]$.
Consider the mapping
\[
\Delta_\zz=(\phi_{z_1}^{-1}\otimes\cdots\otimes\phi_{z_N}^{-1})
\circ\Delta^{(N-1)}:
U(\g[t])\rightarrow U(\g[t])\otimes \cdots\otimes U(\g[t]),
\]
where $\phi_z$ is the automorphism of $U(\g[t])$ changing $t$ to $t-z$.
Here $\zz$ denotes the set of points $\{z_1,\ldots,z_N\}$.

For $x\in\gc$ and $v^{(i)}\in L_i$, we have
\begin{equation}\label{ANTPOD}
(xv^{(0)})\otimes v^{(1)}\otimes\cdots\otimes v^{(N)}
=v^{(0)}\otimes \Delta_\zz(S(x))(v^{(1)}\otimes\cdots\otimes v^{(N)})
\end{equation}
in the quotient space (\ref{CB}). Here $S$ is the antipode.

Let $A\subset U(\g[t])\otimes\cdots\otimes U(\g[t])$ be the left ideal
which annihilates $v_1\otimes\cdots\otimes v_N$.
The kernel of the map (\ref{map}) contains
$\Delta_\zz^{-1}(\ant(A))\Lkinf$.
In fact, we will show that they are equal (see Proposition~\ref{intr}).

In order to prove this statement, we introduce
the notion of fusion right ideal
$X_1 \fusn \dots \fusn X_N (\zz)=\Delta_\zz^{-1}(B)$
for right ideals $X_1,\ldots,X_N\subset U(\g[t])$ by replacing $\ant(A)$ with
\[
B=\sum_{i=1}^N
U(\g[t])^{\otimes(i-1)}\otimes X_i\otimes U(\g[t])^{\otimes(N-i)}.
\]

It is known that the dimension of the space of conformal blocks is given in
terms of  the multiplication rule in the Verlinde algebra (see Section
\ref{VER} for the definition of the Verlinde algebra). In Theorem \ref{vermul},
we give a similar multiplicative formula
for the dimension of the space of coinvariants
\begin{equation}\label{COINV}
\Lkinf/(X_1\fusn\dots\fusn X_N(\zz))\Lkinf
\end{equation}
under a certain condition on the right ideals $X_i$.
In particular, this gives a proof of the equality of
the kernel of (\ref{map}) with the subspace
$(X_1\fusn\dots\fusn X_N(\zz))\Lkinf$
by taking $X_i=\ant(\gX_{\lambda_i}')+B_1$, 
where the left ideal $\gX_{\lambda_i}'\subset U(\g)$
is defined by $U(\g)/\gX_{\lambda_i}'\cong\pi_{\lambda_i}$ and
$B_1=\left(\g\otimes t\C[t]\right) U(\g[t])$.
Therefore, the filtration of the conformal coinvariants (\ref{CB})
can be studied in the form (\ref{COINV}) as the coinvariants with respect to
the fusion right ideal.

Let us consider another
version of this construction. By $\pi(L_i)$ denote the
finite-dimensional representation of $\g$ generated by the highest
weight vector $v_i$. Then one can easily deduce from (\ref{map}) that
the following map is also surjective (see Subsection \ref{subsec:2.4})
\bea\label{map2}
v_0\otimes \pi(L_1) \otimes \dots \otimes \pi(L_N) \to  \left< L_0,
  \dots, L_N\right>^C_{p_0, \dots, p_N}.
\ena

Let $I^{(k)}_\lambda \subset \g \otimes \C[t]$ be the left ideal
which annihilates $v_0$. Then a similar argument shows that 
we have a surjective map
\be
\pi(L_1) \otimes \dots \otimes \pi(L_N) / S(I^{(k)}_\lambda) \to \left< L_0,
  \dots, L_N\right>^C_{p_0, \dots, p_N}.
\en
We will show also that these spaces are in fact isomorphic (see Example~1 
in Subsection~3.2).

To describe the filtration on the space of conformal coinvariants
in this setting, 
we use the notion of filtered tensor product 
introduced in \cite{FL}. It is the tensor product of cyclic
$\g$-modules $V_1, \cdots, V_N$
considered as the evaluation $\g[t]$-modules with distinct parameters $z_i$.
Upon choosing cyclic vectors $v_i \in V_i$ we can 
define a $\g$-equivariant filtration. Namely, the vector $v_1 \otimes
\cdots \otimes v_N$ is cyclic in the tensor product with respect to the
action of $\g[t]$, so the whole representation is a quotient of $U(\g[t])$.
We have a grading of $\g[t]$ and therefore of   
$U(\g[t])$ by the degree in
$t$, so the quotient inherits the corresponding filtration. 

In our case we should choose the highest vectors of $\pi(L_i)$ as
cyclic to obtain the same filtration as given by (\ref{COINV}).
The associated graded
space of the filtered tensor product \cite{FL}
\[
V_1*\cdots*V_N(\zz)=\gr\mathcal{F}_{\zz}(V_1,\cdots,V_N)
\]
is called in \cite{FL} the fusion product. 

For $\g =\slt$ we prove that the Hilbert polynomial of this graded space
is independent of the choice of points and coincides with the level-restricted
Kostka polynomial (for definition see \cite{SS}). This calculation uses the
fusion product. The result of \cite{FF} enables us
to give an upper estimate of the Hilbert polynomial for the quotient of
the filtered tensor product. The resulting formula coincides with the known
fermionic formula for the level-restricted Kostka polynomial. Since at $q=1$,
the level-restricted Kostka polynomial gives the dimension of the conformal
blocks, which is given by the fusion rule, we can conclude that 
the fermionic formula is exact.

The plan of the paper is as follows. In Section 2, we give the definition
of the fusion right ideals in $U(\g[t])$, and prove the fusion
rule for the dimensions of the spaces of coinvariants of integrable
highest weight
$\wg$-modules. In Section 3, we establish an equivalence between spaces of
coinvariants of irreducible highest weight $\wg$-modules and those of filtered
tensor products of cyclic $\g$-modules. We give a brief review on the Kostka
polynomials, and propose the definition of fusion Kostka polynomials.
In Section 4, we prove that in the case of $\g=\slth$,
the associated graded space of the space
of conformal coinvariants is independent of the choice of insertion points, and
that its Hilbert polynomial is the level-restricted Kostka polynomial.
In Appendix we discuss behavior of fusion ideals and spaces of coinvariants when points collide.

\setcounter{section}{1}
\setcounter{equation}{0}

\section{Formula for the dimension}\label{SEC2}
In this section we introduce the notion of fusion right 
ideals and consider spaces of coinvariants with respect to them. 
We give a formula for the dimension of the latter
in terms of the Verlinde algebra (see Theorem \ref{vermul}).

\subsection{Fusion of ideals}
Let $\g$ be a complex Lie algebra. 
By $\gc$ denote the Lie algebra $\g \otimes \C[t]$.
When there is no fear of confusion, we will 
identify $\g$ with $\g\otimes t^0\subset\g[t]$.
For $z\in\C$ we denote by $\phi_z$ the automorphism of $\gc$ sending
$x\otimes t^n$ into $x\otimes (t-z)^n$ $(x\in\g)$. 

Let $\zz = (z_1, \dots, z_N)$ be a set of complex numbers.
Throughout this paper we assume that
\bea
z_i\neq z_j\quad (i\neq j).
\label{distinct}
\ena
Consider the map 
$\g[t] \rightarrow \oplus^N\g[t]$ 
which sends $x\in \g[t]$ to 
$\phi^{-1}_{z_1}(x) \oplus \dots \oplus \phi^{-1}_{z_N}(x)$.
By the universal property it induces the map 
$U(\g[t])\rightarrow U(\g[t])^{\otimes N}$ given by 
\begin{eqnarray*}
\Delta_{\mathcal Z}=\left(\phi_{z_1}\otimes\cdots\otimes\phi_{z_N}\right)^{-1}
\circ\Delta^{(N-1)},
\end{eqnarray*}
where the iterated coproduct is defined by 
$\Delta^{(n-1)}=\left(\Delta\otimes{\rm id}\right)\circ\Delta^{(n-2)}$,  
$\Delta^{(1)}=\Delta$.  
\begin{defn}
{\it Let $X_1, \dots, X_N$ be right ideals in $U(\g[t])$. 
Define the {\em fusion right ideal} by 
$$
X_1 \fusn \dots \fusn X_N (\zz) = \left(\Delta_\zz\right)^{-1}\,
\left(
\sum_{i=1}^NU(\g[t])^{\otimes(i-1)}\otimes X_i\otimes U(\g[t])^{\otimes (N-i)}
\right).
$$}
Throughout this paper, for simplicity of notation, we denote as
\[
I^{(N)}(X_1,\ldots,X_N)=
\sum_{i=1}^NU(\g[t])^{\otimes(i-1)}\otimes X_i\otimes U(\g[t])^{\otimes (N-i)}.
\]
\end{defn}
\medskip

The following Proposition validates this definition.
\begin{prop}\label{prop:2.0}
If the right ideals $X_i=\ga_iU(\g[t])$ are generated by Lie subalgebras
$\ga_i\subset U(\g[t])$, then the fusion right ideal is given by 
\be
X_1 \fusn \dots \fusn X_N (\zz) = 
\Bigl(\bigcap_{i=1}^N\phi_{z_i}(\ga_i)\Bigr)U(\g[t]).
\en
\end{prop}
As an example, consider the two-sided ideal of $U(\g[t])$ 
consisting of elements which have $M$-fold zeros at the origin,  
\bea
B_M = \left(\g \otimes t^M  \C[t]\right) U(\g[t]). 
\label{BM}
\ena
Let 
\bea
B_{M,\zz} = \left(\g \otimes (t-z_1)^M \dots (t-z_N)^M  \C[t]\right) U(\g[t]) 
\label{BM2}
\ena
be the two-sided ideal consisting of elements which have $M$-fold 
zeros at the points in $\zz$. 
Then we have
\be
B_M \fusn \dots \fusn B_M (\zz) = B_{M,\zz}.
\en

For the proof of Proposition \ref{prop:2.0} we use the following fact, 
which is a corollary of the Poincar\'e-Birkhoff-Witt (PBW) theorem. 
\begin{lem}\label{LEMONE}
Let ${\mathfrak a}$ be a Lie algebra, and let
${\mathfrak b}\subset{\mathfrak a}$ be its Lie subalgebra. Take a basis of
${\mathfrak b}$, $\{B_i\}_{i\in I}$, and extend it to a basis of
${\mathfrak a}$ by adjoining $\{A_j\}_{j\in J}$.
We assume that the index set $I$ and $J$ are ordered. 
Then the set of elements in $U({\mathfrak a})$
\begin{equation}\label{SETOF}
B_{i_1}\cdots B_{i_m}A_{j_1}\cdots A_{j_n}
\quad(
i_1\geq\cdots\geq i_m;
j_1\geq\cdots\geq j_n)
\end{equation}
with $m\ge 1,n\ge0$ forms a basis of ${\mathfrak b}U({\mathfrak a})$, 
and the one with $m=0,n\ge 0$ forms a basis of $U(\ga)/\mathfrak{b}U(\ga)$. 
\end{lem}

\begin{lem}\label{LEMTWO}
Let $V_i$ be a vector space with countable basis, 
and let $W_i$ be a subspace of $V_i$, $i=1,\cdots,N$. 
Then the kernel of the map
\[
V_1\otimes\cdots\otimes V_N\rightarrow
(V_1/W_1)\otimes\cdots\otimes(V_N/W_N)
\]
is equal to
\[
\sum_{i=1}^NV_1\otimes\cdots\otimes V_{i-1}\otimes W_i\otimes
V_{i+1}\otimes\cdots\otimes V_N.
\]
\end{lem}
\begin{proof}
Choosing a filtration 
\[
V_i=\bigcup_{a=0}^\infty V_{i,a},\quad V_{i,a}\subset V_{i+1,a}
\]
such that ${\rm dim}\,V_{i,a}$ is finite,
one can reduce the proof to the case where all $V_i$ 
are finite-dimensional. 
Then the statement follows by dimension counting.
\end{proof}
Proposition \ref{prop:2.0} is a consequence
of the following fact. 
\begin{prop}\label{SUBALG}
Let ${\mathfrak a}$ be a Lie algebra with a countable basis, and let
${\mathfrak a}_1,\ldots,{\mathfrak a}_N$ be its Lie subalgebras. Then,
\[
\left(\Delta^{(N-1)}\right)^{-1}\,
\left(
\sum_{i=1}^N U^{\otimes (i-1)}\otimes \ga_i U
\otimes U^{\otimes(N-i)}\right)=(\ga_1\cap\dots\cap\ga_N)U,
\]
where $U=U({\mathfrak a})$.
\end{prop}
\begin{proof}
The inclusion $\supset$ is clear.
Therefore, by using Lemma \ref{LEMTWO}, it is enough to show that the mapping
\[
p:U/(\ga_1\cap\dots\cap\ga_N)U\rightarrow
(U/\ga_1U)\otimes\cdots\otimes(U/\ga_NU),
\]
which is induced from $\Delta^{(N-1)}$, is injective.
We prove this statement by induction on $N$. The case $N=1$ is obvious.
We prove the case $N=2$. Take a basis
$\{A_\alpha\}\sqcup\{B_\beta\}\sqcup\{C_\gamma\}\sqcup\{D_\delta\}$ of $\ga$
such that $\{A_a\}$ is a basis of $\ga_1\cap\ga_2$,
$\{A_a\}\sqcup\{B_b\}$ is a basis of $\ga_1$,
and $\{A_a\}\sqcup\{C_c\}$ is a basis of $\ga_2$.
We fix an ordering of each set
$\{A_\alpha\},\{B_\beta\},\{C_\gamma\},\{D_\delta\}$.
Write $x\in U$ as a linear combination of the PBW basis using the ordering:
\begin{equation}\label{PBWEXP}
x=\sum
c_{\alpha_1,\ldots,\delta_d}
A_{\alpha_1}\cdots A_{\alpha_a}
B_{\beta_1}\cdots B_{\beta_b}
C_{\gamma_1}\cdots C_{\gamma_c}
D_{\delta_1}\cdots D_{\delta_d}
\end{equation}
Here and after, we do not bother to mention a necessary ordering
for the indices.

Let
\[
\pi:U\otimes U\rightarrow(U/\ga_1U)\otimes(U/\ga_2U)
\]
be the canonical surjection. We want to show that if
$\pi\circ\Delta(x)=0$ then $x$ belongs to $(\ga_1\cap\ga_2)U$.

Let $F^lU$ be the filtration of $U$ defined as follows.
\[
F^0U=\C\cdot1,\quad F^{l+1}U=F^lU+{\mathfrak a}F^lU.
\]
One can define a filtration of $U\otimes U$ by setting
\[
F^l(U\otimes U)=\sum_{l_1+l_2=l}F^{l_1}U\otimes F^{l_2}U,
\]
and further induce a filtration of $(U/\ga_1U)\otimes(U/\ga_2U)$.
By Lemma \ref{LEMONE}, the elements of the form
$C_{\gamma_1}\cdots C_{\gamma_c}D_{\delta_{i_1}}\cdots D_{\delta_{i_{d'}}}$
(resp., $B_{\beta_1}\cdots B_{\beta_b}D_{\delta_{j_1}}\cdots D_{\delta_{j_{d''}}}$)
form a basis of $U/\ga_1U$ (resp., $U/\ga_2U$).
We have the basis of $(U/\ga_1U)\otimes(U/\ga_2U)$ consisting of
the tensor product of these bases. Moreover, those which satisfies
$c+b+d'+d''=l$, form a basis of
$(F^l/F^{l-1})((U/\ga_1U)\otimes(U/\ga_2U))$.
Let us denote them by 
$X_{\gamma_1,\ldots,\gamma_c,\delta_{i_1}.\ldots,\delta_{i_{d'}};
\beta_1,\ldots,\beta_b,\delta_{j_1},\ldots,\delta_{j_{d''}}}$.

Note that $\pi\circ\Delta(1)\not=0$. Suppose that
$x\in F^lU$ $(l\geq1)$ and $\pi\circ\Delta(x)=0$.
We will show that $x\in(\ga_1\cap\ga_2)U+F^{l-1}U$.
Write $x$ as (\ref{PBWEXP}). If a monomial in this expression is
such that $a\geq1$ then it belongs to $(\ga_1\cap\ga_2)U$ and
the kernel of $\pi\circ\Delta$. Therefore,
one can assume that only the terms with $a=0$ appear in (\ref{PBWEXP}).

We want to show that
\begin{equation}\label{EN1}
c_{\beta_1,\ldots,\delta_d}=0
\end{equation}
if $b+c+d=l$.
Consider the $2^l$ monomials that are obtained by expanding the product
$\Delta(B_{\beta_1})\cdots\Delta(B_{\beta_b})
\Delta(C_{\gamma_1})\cdots\Delta(C_{\gamma_c})
\Delta(D_{\delta_1})\cdots\Delta(D_{\delta_d})$.
If a monomial is such that either one of $B_\beta$ is sent to the first
component of the tensor product $U\otimes U$, or one of $C_\gamma$
to the second, then the image of that monomial by the surjection $\pi$
belongs to $F^{l-1}((U/\ga_1U)\otimes(U/\ga_2U))$. Therefore,
$\pi\circ\Delta(B_{\beta_1}\cdots D_{\delta_d})\in
(F^l/F^{l-1})((U/\ga_1U)\otimes(U/\ga_2U))$
is the sum of $2^d$ elements
$X_{\gamma_1,\ldots,\gamma_c,\delta_{i_1}.\ldots,\delta_{i_{d'}};
\beta_1,\ldots,\beta_b,\delta_{j_1},\ldots,\delta_{j_{d''}}}$
where
$\{i_1.\ldots,i_{d'}\}\sqcup\{j_1.\ldots,j_{d''}\}$
is a partition of the index set $\{1,\ldots,d\}$.
Let us denote this sum by $Y_{\beta_1,\ldots,\delta_d}$.
The elements $Y_{\beta_1,\ldots,\delta_d}$ that are obtained
from $x$ are linearly independent. Therefore, we have (\ref{EN1}).

Finally, suppose that the statement is true for $N-1$ $(N\geq3)$. It means
that the mapping
\[
U/(\ga_1\cap\dots\cap\ga_{N-1})U\rightarrow
(U/\ga_1U)\otimes\cdots\otimes(U/\ga_{N-1}U)
\]
is injective. From this follows that
\[
p_1:(U/(\ga_1\cap\dots\cap\ga_{N-1})U)\otimes(U/\ga_NU)\rightarrow
(U/\ga_1U)\otimes\cdots\otimes(U/\ga_NU)
\]
is also injective. Consider the mapping
\[
p_2:U/(\ga_1\cap\dots\cap\ga_N)U\rightarrow
(U/(\ga_1\cap\dots\cap\ga_{N-1})U)\otimes(U/\ga_NU)
\]
induced from $\Delta$. Since
$\Delta^{(N-1)}=(\Delta^{(N-2)}\otimes{\rm id})\circ\Delta$,
we have $p=p_1\circ p_2$. Therefore, it is enough to show that
$p_2$ is injective. This is nothing but the case $N=2$.
\end{proof}

Note that $\gc$ and therefore $U(\gc)$ are graded by the degree in $t$.

\begin{defn}
{\it 
We call a right ideal $X \subset U(\g[t])$ a {\em congruence ideal} if X is 
graded and $X \supset B_M$ for a sufficiently large $M$.
}
\end{defn}

Quite generally, let $V$ be a left module over a Lie algebra $\ga$.
For a left ideal $X\subset U(\ga)$, we denote by
$V^X=\{v\in V \mid xv=0~~\forall x\in X\}$ the space of $X$-invariants. 
Likewise, for a right ideal $Y \subset U(\ga)$, we denote by $V/Y=V/YV$
the space of $Y$-coinvariants. We have a canonical isomorphism
\bea
(V/Y)^*\simeq (V^*)^{S(Y)},
\label{iso1}
\ena
where $V^*=\mathop{\rm Hom}_\C(V,\C)$ and $\ant:U(\ga)\to U(\ga)$ stands for
the antipode (the anti-automorphism such that $\ant(x)=-x$ for $x\in\ga$). 
We will consider coinvariant spaces with respect to fusion right ideals. 

\subsection{Affine Lie algebras}
In this subsection we fix our notation concerning affine Lie algebras and integrable modules. 

Let $\g$ be a simple Lie algebra. 
Fix a triangular decomposition $\g=\gn_+\oplus\hh\oplus\gn_-$, 
where $\hh \subset \g$ is the Cartan subalgebra  and 
$\gn_+$ (resp. $\gn_-$) is the nilpotent subalgebra of all positive (resp. negative) root vectors. 
The invariant bilinear form $(~|~)$ on $\g$ is normalized as 
$(\theta|\theta)=2$, where $\theta$ is the maximal root. 
Let $P$ be the weight lattice and $P_+$ be the set of dominant integral weights.  
For a non-negative integer $k\in\Z_{\ge 0}$, we set 
$\simk=\{\lambda\in P_+\mid \br{\lambda,\theta^{\vee}}\le k\}$,
where $\theta^{\vee}$ denotes the maximal coroot. 
For $\lambda \in P_+$, 
we denote by $\pi_\lambda$ the irreducible representation with 
highest weight $\lambda$.
We will write $\lambda^*=-w_0(\lambda)$ for 
the highest weight of the dual representation $\pi_\lambda^*$,
$w_0$ being the longest element of the Weyl group $W$. 
If $\lambda\in P^{(k)}_+$, then $\lambda^*\in P_+^{(k)}$.

Let $\wg = \g\otimes\C[t,t^{-1}]\oplus \C K\oplus\C d$ be 
the affine Lie algebra of non-twisted type.  
Here $K$ is a central element, $d=td/dt$, 
\be
&&[x\otimes f,y\otimes g]=[x,y]\otimes fg
+\Omega_0(x\otimes f,y\otimes g)K 
\\
&&\qquad
\qquad\qquad \qquad (x,y\in\g, f,g\in\C[t,t^{-1}]),
\en
and the canonical cocycle is given by
\be
\Omega_0(x\otimes f,y\otimes g)=(x|y)\mathop{\rm Res}_{t=0}df\cdot g.
\en
Denote by $\Lambda_i$ the fundamental weights
$(0\le i\le r;r=\dim\mathfrak{h})$ and by $\delta$ the null root. 
Let $\Lambda=\lambda+k\Lambda_0+a\delta$ 
$(\lambda\in P$, $a\in\C)$ be
an affine weight of level $k$.
For an element $w$ of the affine Weyl group $W_{\rm aff}$, 
we set 
\bea
&&w\circ(\lambda+k\Lambda_0+a\delta)=w(\lambda,k)+k\Lambda_0
+(a-d_w(\lambda,k))\delta,
\label{shift}
\\
&&w(\lambda,k)\in P, d_w(\lambda,k)\in\Z, 
\nn
\ena
where 
$w\circ \Lambda=w(\Lambda+\rho)-\rho$ 
stands for the shifted action
and $\rho=\sum_{i=0}^r\Lambda_i$. 

We say a $\widehat{\g}$-module has level 
$k$ if $K$ acts as $k$ times the identity on it. 
For $k\in\Z_{\ge 0}$ and $\lambda \in P_+^{(k)}$, 
let $\Lkla=U(\widehat{\g})\vkla$ be the unique level $k$ 
integrable module such that 
\bea
&&(\g\otimes t^n)\vkla=0 \quad(n>0),
\label{high1}\\
&&\gn_+ \vkla=0,
\label{high2}\\
&&
h\vkla=\br{\lambda,h}\vkla\quad (h\in\mathfrak{h}),
\quad 
d\vkla=0.
\label{high3}
\ena

We also consider modules `placed at infinity'. 
To be precise, define 
$\wg(\infty) = \g\otimes\C[t,t^{-1}]\oplus \C K\oplus\C d$ by 
the same equations as above, except that we change the cocycle to 
\be
\Omega_\infty(x\otimes f,y\otimes g)
&=&(x|y)\mathop{\rm Res}_{t=\infty}df\cdot g
\\
&=&-\Omega_0(x\otimes f,y\otimes g).
\en
Let $\Lkinf=U(\wg(\infty))\vkinf$ be the integrable 
$\wg(\infty)$-module characterized by \eqref{high1}--\eqref{high3}, 
with $\vkinf$ in place of $\vkla$ and 
\eqref{high1} being changed to
\be
(\g\otimes t^n)\vkinf=0 \quad (n<0).
\en
The degree zero part of $\Lkla$ and $\Lkinf$ are both isomorphic to
$\pi_\lambda$ as $\g$-module.

Consider the isomorphism of Lie algebra $\iota:\wg\rightarrow\wg(\infty)$ 
given by $x\otimes t^n\mapsto x\otimes t^n$, $K\mapsto -K$, $d\mapsto d$.
By this isomorphism, $\Lkinf$ can be considered as the level $-k$
lowest weight module with the lowest weight $-\lambda$.
We have a non-degenerate pairing 
\bea
\Lkinf\times L^{(k)}_{\lambda^*}\longrightarrow \C
\label{pair}
\ena
such that $\br{\iota(x)u,v}=-\br{u,xv}$ 
$(x\in\wg,u\in \Lkinf,v\in L^{(k)}_{\lambda^*})$, which extends
the canonical coupling $\pi_\lambda\times \pi_{\lambda^*}\rightarrow\C$.

\begin{prop}\label{prop:2.1}
Let $X\subset U(\g[t])$ be a congruence right ideal. 
Then $\Lkinf/X$ is finite-dimensional. We have an isomorphism 
\bea
(\Lkinf/X)^* \cong (L^{(k)}_{\lambda^*})^{\ant(X)}.
\label{iso0}
\ena
\end{prop}

\begin{proof}
First we prove the finite dimensionality of $\Lkinf/X$. 
Choose a set of nilpotent generators $x_1,\cdots,x_m$ of $\g$, 
and let $\widetilde{B}_{M,i}$ be the right ideal generated by
$x_j\otimes t^n$ $(1\le j\le m, n\ge M)$ and 
$x_1\otimes t^{M-1},\cdots,x_i\otimes t^{M-1}$.
Note that $\widetilde{B}_{M,0}=\widetilde{B}_{M+1,m}$. We have
$\widetilde{B}_{M,i}=\widetilde{B}_{M,i-1}+(x_i\otimes t^{M-1})U(\g[t])$.
Since the congruence ideal $X$ contains $\widetilde{B}_{M,0}$
for a large $M$, in order to show the finite-dimensionality of
$\Lkinf/X$, it suffices to prove it for  
$\Lkinf/\widetilde{B}_{M,i}$ for all $M,i$.
Since $\Lkinf/\widetilde{B}_{1,0}\simeq \pi_\lambda$,
the statement is valid for $M=1,i=0$. The general case
can be shown by induction on $M,i$, by applying the following lemma
to $V=\Lkinf/\widetilde{B}_{M,i-1}$, $x=x_i\otimes t^{M-1}$.

Eq.\eqref{pair} states that $L^{(k)}_{\lambda^*}$ 
is the restricted dual of $\Lkinf$. Since $X$ is graded
and $\Lkinf/X$ is finite-dimensional, \eqref{iso0} follows from \eqref{iso1}.
\end{proof}

\begin{lem}\label{lem:2.1} 
Let $V$ be a vector space and $x\in\End(V)$.
Assume that $x$ is locally 
nilpotent, $\bigcap_{n>0}\mathop{\rm Im}x^n=\{0\}$, and
that $V/\mathop{\rm Im}x$ is finite dimensional. 
Then $V$ is also finite dimensional.
\end{lem}

\begin{proof}
Choose vectors $v_1,\cdots,v_n\in V$ which span $V/\mathop{\rm Im}x$, and let $M$ 
be such that $x^Mv_i=0$ for all $i$. 
Then the set $\{x^jv_i\}$ with $0\le j\le M-1,1\le i\le n$ spans $V$.
\end{proof}

\subsection{The multiplicative formula for dimensions}\label{VER}
The aim of this section is to prove the multiplicative formula
for the dimension of the space of coinvariants with respect to the
fusion right ideals. Before discussing it, we recall the original version of
the fusion rule for the space of conformal coinvariants.

{}From now on, we fix a level $k\in\Z_{\ge 0}$. 
First recall the definition of the Verlinde algebra.
Let $\mathop{\rm Rep}\g=\oplus_{\lambda\in P_+}\Z\cdot[\lambda]$ be 
the Grothendieck ring of $\g$, where $[\lambda]$ stands
for the class of irreducible representation $\pi_\lambda$. 
Let $\mathcal{R}^{(k)}$ denote the ideal spanned by elements of the form
\be
[\lambda]-(-1)^{l(w)}[w(\lambda,k)]
\qquad 
(w\in W_{\rm aff},\lambda,w(\lambda,k)\in P_+),
\en
where $l(w)$ signifies the length of $w$. The Verlinde algebra is
the quotient ring $\mathcal{V}^{(k)}=\mathop{\rm Rep}\g/\mathcal{R}^{(k)}$.  
We use the same letter $[\lambda]$ to represent its class
in $\mathcal{V}^{(k)}$.  The set $\{[\lambda]\}_{\lambda\in P_+^{(k)}}$
constitutes a $\Z$-basis of $\mathcal{V}^{(k)}$. For an element
$a\in\mathcal{V}^{(k)}$, we write the coefficient of $[\lambda]$ 
in this basis as $(a:[\lambda])_k$, i.e., 
\be
a=\sum_{\lambda\in P_+^{(k)}}(a:[\lambda])_k[\lambda].
\en

We now fix a set of distinct points $\zz=\{z_1,\cdots,z_N\}$ in $\C$.
\begin{thm}\cite{TUY}\label{thm:verlinde}
Let $L_i$ $(i=0,\ldots,N)$ be the irreducible highest weight representation
of highest weight $\lambda_i$. Define the space of conformal coinvariants
by (\ref{CB}). Then, we have the dimension formula.
\begin{equation}\label{VERLINDE}
{\rm dim}\,
\left<L_0,\dots,L_N\right>^{\C P^1}_{\infty,z_1,\dots,z_N}=
([\lambda_1^*]\cdots[\lambda_N^*]:[\lambda_0])_k.
\end{equation}
\end{thm}

Now, consider a set of congruence right ideals $X_1,\cdots,X_N$ of
$U(\g[t])$. 
\begin{lem}
The space of coinvariants
\bea
\Lkinf/X_1\fusn\cdots\fusn X_N(\zz)
\label{coinv2}
\ena
is finite dimensional.
\end{lem}
\begin{proof}
A standard deformation argument (e.g.\cite{FKLMM1}, Lemma 20)
shows that 
$\dim \Lkinf/B_{M,\zz}\le \dim \Lkinf/B_{NM}<\infty$, 
where we used Proposition \ref{prop:2.1}. 
Since $X_1\fusn\cdots\fusn X_N(\zz)\supset B_{M,\zz}$ for some $M$, 
we obtain the assertion. 
\end{proof}
The main result of this section is the following. 
\begin{thm}\label{vermul}
Notation being as above, we have
\be
&&\sum_{\lambda \in \simk} 
\dim \bigl(\Lkinf/X_1\fusn\dots\fusn X_N(\zz)\bigr)\cdot[\lambda]
\\
&& 
\qquad =  
\prod_{i=1}^N
\Bigl(\sum_{\lambda \in \simk} \dim \bigl(\Lkinf/X_i\bigr) 
\cdot [\lambda]\Bigr).
\en
\end{thm}

For the proof, we make use of a construction due to \cite{FFu}
which we recall below.
See also appendix to \cite{FKLMM1}, where
a detailed account is given in the case $\g=\slt$. 

Set 
$$ 
\Cc^{(k)}_{\lambda,M}(\zz) = \Lkinf /B_{M,\zz} 
$$
where $B_{M,\zz}$ is defined in \eqref{BM2}.
Since $B_{M,\zz}$ is a two-sided ideal, 
$\Cc^{(k)}_{\lambda,M}(\zz)$ is a $U(\g[t])$ module.  
Taking the dual module, we obtain 
an inductive system
$\cdots\subset (\Cc^{(k)}_{\lambda,M}(\zz))^*\subset
(\Cc^{(k)}_{\lambda,M+1}(\zz))^*\subset\cdots$.
The direct limit
\be
\Lc^{(k)}_{\lambda}(\zz) = \varinjlim_M\left(\Cc^{(k)}_{\lambda,M}(\zz)^*\right)
\en
admits an action of the algebra 
$$
U_\zz(\gc) = 
\varprojlim_M U(\gc)/B_{M,\zz}.
$$
Note that, in the special case where $\zz=\{0\}$, 
\be
U_{\{0\}}(\g[t])=
\varprojlim_M U(\gc)/\left(t^M\g[t]U(\g[t])\right)
\en
is the completion by formal power series in $t$.
We denote the completion of $U(\gc)^{\otimes N}$ as
\begin{eqnarray}\label{COMPL}
&&\qquad U_{\{0\}}(\gc)^{\widehat\otimes N}\\
&&\qquad=\varprojlim_M
U(\gc)^{\otimes N}/I^{(N)}(B_M,\ldots,B_M).\nn
\end{eqnarray}

We use the following lemmas.
\begin{lem}\label{LEM2}
Suppose ${\mathfrak b}$ is an ideal of a Lie algebra ${\mathfrak a}$.
Then, we have
\[
U({\mathfrak a}/{\mathfrak b})\cong
U({\mathfrak a})/{\mathfrak b}U({\mathfrak a}).
\]
\end{lem}
\begin{lem}\label{LEM3}
Suppose ${\mathfrak a}={\mathfrak a}_1\oplus\cdots\oplus{\mathfrak a}_N$
is a direct sum decomposition of a Lie algebra.
Then, we have an isomorphism of algebras
\[
U({\mathfrak a})\cong
U({\mathfrak a}_1)\otimes\cdots\otimes U({\mathfrak a}_N).
\]
\end{lem}

Now, we state a basic proposition about localization of
the universal enveloping algebra.
\begin{prop}\label{loc}
$(i)$
We have a canonical isomorphism 
\be
U_\zz(\gc) \cong \cpu^{\widehat\otimes N}. 
\en

$(ii)$
The canonical inclusion
$U=U(\gc)\to U_\zz(\gc)\cong\cpu^{\widehat\otimes N}$ is the composition of 
$\Delta_\zz$ followed by the inclusion 
$U(\gc)^{\otimes N}\to\cpu^{\widehat\otimes N}$.

$(iii)$
The action of $U^{\otimes N}$ on $\Lc^{(k)}_{\lambda}(\zz)$
can be extended to an action of $U(\wg)^{\otimes N}$.
\end{prop}
\begin{proof}
Since the points in $\zz$ are distinct, the restriction of 
$\Delta_\zz$ to $\g[t]$ induces an isomorphism
\be
\iota_\zz:
\g[t]/\prod_{i=1}^N(t-z_i)^M\g[t]
\cong
\oplus^N\left(\g[t]/t^M\g[t]\right).
\en
By Lemma \ref{LEM3} we have
\[
\Delta_\zz:U\left(\gc/\prod_{i=1}^N(t-z_i)^M\gc\right)
{\buildrel\cong\over\longrightarrow}U\!\left(\gc/t^M\gc\right)^{\otimes N}
\]
By Lemma \ref{LEM2}, we have
\[
U/B_{M,\zz}\cong\left(U/t^M\gc U\right)^{\otimes N}.
\]

For (i) and (ii), it is enough to show that the kernel of the canonical
surjection $U^{\otimes N}\rightarrow (U/t^M\gc U)^{\otimes N}$ coincides with
$I^{(N)}(B_M,\ldots,B_M)$.
This follows from Lemma \ref{LEMTWO}. 

The action of $U(\widehat{\g})^{\otimes N}$ is defined as follows. 
Consider an element $x(t)=(x_1(t),\cdots,x_N(t))$ of $\oplus^N\wg$. 
Choose an $m$ such that $t^m x_i(t)\in \g[t]$ for all $i$. 
Then $\prod_{j=1}^N(t-z_j)^m\iota_\zz^{-1}(x(t))\in\g[t]$, so that 
$\iota_\zz^{-1}(x(t))$ maps 
$\Cc^{(k)}_{\lambda,M}(\zz)$  to 
$\Cc^{(k)}_{\lambda,M-m}(\zz)$ for $M\ge m$. 
Therefore $x(t)$ acts on $\Lc^{(k)}_{\lambda}(\zz)$. 
\end{proof} 

We will use the following fact. 
\begin{thm}\label{thm:FFu}\cite{FFu}
We have an isomorphism of $U(\wg)^{\otimes N}$-modules.
$$
\Lc^{(k)}_{\lambda}(\zz) \cong \bigoplus_{\lambda_1, \dots, \lambda_N \in \simk}
([\lambda_1]\cdot \dots \cdot[\lambda_N]: [\lambda])_k \cdot 
L^{(k)}_{\lambda_1^*} 
\boxtimes \dots \boxtimes L^{(k)}_{\lambda_N^*}.
$$
\end{thm}

For $I \subset U(\gc)$, denote  by $I^{\zz}$
its closure in $U_\zz(\gc)$ in the topology of the inverse limit.

\begin{lem}\label{lem:2.3} 
$$
\left(\Lkinf/ X_1 \fusn \dots \fusn X_N(\zz)\right)^* \cong 
\left(\Lc^{(k)}_{\lambda}(\zz)\right)^{\ant(X_1\fusn \dots \fusn X_N(\zz))^\zz}.
$$
\end{lem}

\begin{proof}
Choose $M$ such that $X_i \supset B_M$ for all $i$. 
Then
$$
\Lkinf/ X_1 \fusn \dots \fusn X_N(\zz) 
\cong
\Cc^{(k)}_{\lambda,M}(\zz) / X_1 \fusn \dots \fusn X_N(\zz).
$$
We have
$$
\left(\Cc^{(k)}_{\lambda,M}(\zz)/X_1\fusn\dots\fusn X_N(\zz)\right)^*\cong
\left(\Cc^{(k)}_{\lambda,M}(\zz)^*\right)^{\ant(X_1\fusn\dots\fusn X_N(\zz))}.
$$
Taking the direct limit in $M$ we obtain the statement of the lemma. 
\end{proof}

\begin{lem}\label{LEM4}
Suppose that $z_1,\ldots,z_N\in\C$ are distinct.
For any $L>0$ and $i=1,\cdots,N$, there exists an $f_i(t)\in\C[t]$ satisfying
\[
f_i(t)\equiv
\begin{cases}
1\bmod(t-z_i)^L;\\
0\bmod(t-z_j)^L&\hbox{for }j\not=i.\\
\end{cases}
\]
\end{lem}

\begin{lem}\label{lem:2.4}  
Under the identification of Proposition~\ref{loc} we have
$$
(X_1\fusn \dots \fusn X_N(\zz))^\zz \cong 
\sum_{i=1}^NU_{\{0\}}(\gc)^{\widehat\otimes (i-1)}
\widehat\otimes X_i^{\{0\}}\widehat\otimes
U_{\{0\}}(\gc)^{\widehat\otimes (N-i)}.
$$
Here, the right hand side is the completion in (\ref{COMPL}).
\end{lem}
\begin{proof}
For a sufficiently large $L$, we have 
\begin{equation}\label{NOTLOK}
I^{(N)}(B_L,\ldots,B_L)\subset I^{(N)}(X_1,\ldots,X_N).
\end{equation}
Recall that
$X_1\fusn \dots \fusn X_N(\zz)=\Delta_\zz^{-1}(I^{(N)}(X_1,\ldots,X_N))$.

We want to show that for a sufficiently large $L$ and for any element
$a\in I^{(N)}(X_1,\ldots,X_N)$,
there exists an element $x\in X_1\fusn \dots \fusn X_N(\zz)$ satisfying
\begin{equation}\label{WANT}
\Delta_\zz(x)\equiv a\bmod I^{(N)}(B_L,\ldots,B_L).
\end{equation}
{}From (\ref{NOTLOK}) it is enough to find
$x\in U(\gc)$ satisfying (\ref{WANT}) because it implies
$x\in X_1\fusn \dots \fusn X_N(\zz)$. The existence of such an 
$x$ follows from 
\[
\Delta_\zz(f_i(t)g)\equiv1\otimes\cdots\otimes
\overset{\scriptstyle{i-{\rm th}}}{g}
\otimes\cdots\otimes1\bmod I^{(N)}(B_L,\ldots,B_L),
\]
where $g\in\gc$ and $f_i(t)$ is given in Lemma \ref{LEM4}.
\end{proof}

{\it Proof of Theorem \ref{vermul}.}
Set $X=X_1\fusn \dots \fusn X_N(\zz)$. 
Collecting together the statements above, we find
\begin{eqnarray*}
(\Lkinf/ X)^* 
&\cong &
(\Lc^{(k)}_{\lambda}(\zz))^{\ant(X)^\zz} 
\qquad\qquad\qquad 
\mbox{(by Lemma \ref{lem:2.3})}
\\
&
\cong&\hskip-15pt \bigoplus_{\lambda_1, \cdots, \lambda_N \in \simk} 
\hskip-10pt([\lambda_1] \dots [\lambda_N]: [\lambda])_k \cdot
(L^{(k)}_{\lambda_1^*} \boxtimes \dots \boxtimes L^{(k)}_{\lambda_N^*})^{
\ant(X)^\zz
}
\\
&&\qquad\qquad\qquad\qquad\qquad \qquad 
\mbox{(by Theorem \ref{thm:FFu})}
\\
&
\cong&\hskip-20pt \bigoplus_{\lambda_1, \cdots, \lambda_N \in \simk} 
\hskip-10pt([\lambda_1] \dots [\lambda_N]: [\lambda])_k \cdot
(L^{(k)}_{\lambda_1^*})^{\ant(X_1)} \boxtimes \dots \boxtimes
(L^{(k)}_{\lambda_N^*})^{\ant(X_N)} 
\\
&&\qquad\qquad\qquad\qquad\qquad \qquad \mbox{(by Lemma \ref{lem:2.4})}
\\
&
\cong&\hskip-15pt \bigoplus_{\lambda_1, \dots, \lambda_N \in \simk}  
([\lambda_1] \dots [\lambda_N]: [\lambda])_k
\\
&&\qquad\qquad\times
(L^{(k)}_{\lambda_1}(\infty)/X_1 \boxtimes \dots \boxtimes
L^{(k)}_{\lambda_N}(\infty)/X_N)^*
\\
&&\qquad\qquad\qquad\qquad\qquad \qquad \mbox{(by Proposition \ref{prop:2.1}).}
\end{eqnarray*}
Therefore
$$
\dim \Lkinf/ X
= \sum_{\lambda_1,\dots,\lambda_N\in\simk}   
([\lambda_1] \cdots [\lambda_N]: [\lambda])_k \cdot
\prod_{i=1}^N \dim L^{(k)}_{\lambda_i}(\infty) / X_i,
$$
which proves the theorem.
\qed

\subsection{Space of conformal coinvariants}\label{subsec:2.4}
We give a sketch of proof for the surjectivity statements in Introduction,
which motivate the definition of spaces of coinvariants in two different
settings. In Section \ref{sec:3} these two spaces will be identified.
We also conclude Section \ref{SEC2} by stating the isomorphism
between the space of conformal coinvariants and the space of coinvariants
by the fusion right ideal.

The proof of the surjectivity of (\ref{map}) and (\ref{map2}) is based
on the equality
(\ref{ANTPOD}) where we take $x\in\g^{C\setminus\{p_0,\dots,p_N\}}$.
There exists a meromorphic function $f$ on the curve $C$ which has poles only
at $p_0$ and $p_i$ $(i\not=0)$. Moreover, the singular part
of the pole at $p_i$
can be chosen arbitrarily. By using (\ref{ANTPOD}) for $x=X\otimes f$
$(X\in\g)$, one can prove inductively the surjectivity of the map
\[
L_0\otimes \pi(L_1)\otimes\dots\otimes\pi(L_N)\to
\left< L_0,\dots,L_N\right>^C_{p_0,\dots,p_N}.
\]
Now, for (\ref{map}), taking $f$ which has a pole only at $p_0$ and
such that $f(p_j)=\delta_{i,j}$, we can reduce the vectors in $\pi(L_i)$
to the highest weight vector $v_i$. For (\ref{map2}), we assume that the curve
$C$ is rational. Then, taking $f$ which has a pole only at $p_0$
(we can choose an arbitrary singular part of the pole),
we can reduce the vectors in $L_0$ to the highest weight vector $v_0$.

Note that in the above proof we used only the property of the highest weight
vector $v_i$ that it is annihilated by $\g\otimes t_i\C[t_i]$ and
that it is a cyclic vector in $L_i$. Therefore, we can replace $L_i$
by a direct sum of highest weight representations with a cyclic vector $v_i$
in the degree (with respect to $d$) zero component.
The surjectivity (\ref{map}) is still valid.

For $\lambda\in P_+$, we define the left ideal $\gX_\lambda\subset U(\g)$ by
\begin{equation}
\gX_\lambda=
U(\g)\gn_++\sum_{h\in\hh}U(\g)(h-\left<\lambda,h\right>).
\label{gxl}
\end{equation}
We also define $\gX'_{\lambda}\subset U(\g)$ by 
\begin{equation}
\gX_\lambda'=
U(\g)\gn_++\sum_{h\in\hh}U(\g)(h-\left<\lambda,h\right>)+
\sum_{j} U(\g)f_j^{\langle \lambda, \alpha^{\vee}_j \rangle+1}
\label{gxl2}
\end{equation}
where $\alpha_j$ and $\alpha^{\vee}_j$ denote the simple roots and coroots
respectively, and the $f_j$ are root vectors corresponding to $-\alpha_j$. 
Namely $\gX_{\lambda}'$ is the annihilating ideal of the highest weight 
vector $v_\lambda\in\pi_\lambda$. 

\begin{lem}\label{LEMM}
We have
\[
{\rm dim}\,\pi_\mu/S(\gX_\lambda)
={\rm dim}\,\pi_\mu/S(\gX_\lambda') = \delta_{\mu,\lambda^*}.
\]
\end{lem}
\begin{proof} Let $\gX$ be $\gX_\lambda$ or $\gX_\lambda'$.

Then the right ideal $S(\gX)$ contains ${\mathfrak n}_+$
and $h+\langle\lambda,h\rangle$ $(h\in{\mathfrak h})$.
Since $-\lambda=w_0(\lambda^*)$, this implies
$\pi_\mu/S(\gX)=0$ unless $\mu=\lambda^*$.

For $\mu=\lambda^*$, denote by $v^*_\lambda$
the lowest weight vector of $\pi_{\lambda^*}$
such that $\langle v^*_\lambda,v_\lambda\rangle=1$
with the highest weight vector $v_\lambda\in \pi_\lambda$.
We have $\pi_{\lambda^*}/S(\gX)=\C v^*_\lambda$
if one can show that $v^*_\lambda$ is not contained in the image
of $S(\gX)$. Suppose that $v^*_\lambda=S(g)w$ for some
$g\in\gX$ and $w\in \pi_{\lambda^*}$. Then, we have
\[
1=\langle v^*_\lambda,v_\lambda\rangle
=\langle S(g)w,v_\lambda\rangle
=\langle w,gv_\lambda\rangle=0
\]
as elements of $\gX$ annihilate the highest weight vector.
This is a contradiction.
\end{proof}

\begin{cor}\label{CORO}
Set 
$$X_{\lambda_i}=S(\gX_{\lambda_i})+B_1, \qquad
X_{\lambda_i}'=S(\gX_{\lambda_i}')+B_1$$ 
for $\lambda_i\in P^{(k)}_+$.
Then, we have
\[
{\rm dim}\Lkinf/X_{\lambda_1}\fusn\dots\fusn X_{\lambda_N}(\zz)
=\bigl([\lambda_1^*]\cdots[\lambda_N^*]:[\lambda]\bigr)_k,
\]
\[
{\rm dim}\Lkinf/X_{\lambda_1}'\fusn\dots\fusn X_{\lambda_N}'(\zz)
=\bigl([\lambda_1^*]\cdots[\lambda_N^*]:[\lambda]\bigr)_k.
\]
\end{cor}
\begin{proof}
As above, let $\gX$ be $\gX_\lambda$ or $\gX_\lambda'$ and let $X$ be
$X_\lambda$ or $X_{\lambda}'$ respectively.
Then
\[
\Lkinf/X = \Lkinf/(B_1+S(\gX))\cong\pi_\lambda/S(\gX).
\]
The statement follows from Lemma \ref{LEMM} and Theorem \ref{vermul}.
\end{proof}

{}From Corollary \ref{CORO} we have
\begin{prop}\label{intr} 
Let $z_1,\cdots, z_N$ be distinct points in $\C$. 
Set
$$X_{\lambda_i}=S(\gX_{\lambda_i})+B_1, \qquad
X_{\lambda_i}'=S(\gX_{\lambda_i}')+B_1$$
\begin{enumerate}
\item We have the canonical isomorphism given by (\ref{map})
\begin{equation}\label{isoprop}
L^{(k)}_\lambda(\infty)/X_{\lambda_1}'\fusn\cdots\fusn
X_{\lambda_N}'(\zz)
\simeq
\left<L^{(k)}_\lambda,
  L^{(k)}_{\lambda_1},\dots,L^{(k)}_{\lambda_N}\right>^{\C
    P^1}_{\infty, z_1,\dots,z_N}
\end{equation}
\item We also have the canonical isomorphism
\begin{equation}\label{isoprop2}
L^{(k)}_\lambda(\infty)/X_{\lambda_1}\fusn\cdots\fusn
X_{\lambda_N}(\zz)
\simeq
\left<L^{(k)}_\lambda,  
  L^{(k)}_{\lambda_1},\dots,L^{(k)}_{\lambda_N}\right>^{\C
    P^1}_{\infty, z_1,\dots,z_N}
\end{equation}
\end{enumerate}
\end{prop}
\begin{proof} First let us prove (i).
The map (\ref{map}) gives a surjection from 
the left hand side to the right hand side of 
\eqref{isoprop}.
The injectivity follows since their dimensions 
are equal 
by Theorem \ref{vermul} and Corollary \ref{CORO}.

To prove (ii) note that as $X_\lambda \subset X_\lambda'$ we have the
natural
surjection
\begin{equation}\label{surjprop}
L^{(k)}_\lambda(\infty)/X_{\lambda_1}\fusn\cdots\fusn
X_{\lambda_N}(\zz) \to
L^{(k)}_\lambda(\infty)/X_{\lambda_1}'\fusn\cdots\fusn
X_{\lambda_N}'(\zz).
\end{equation}
Theorem \ref{vermul} and Corollary \ref{CORO} imply that the
dimensions of these spaces are equal, so \eqref{surjprop} and therefore
\eqref{isoprop2} is an
isomorphism.
\end{proof}

\setcounter{section}{2}
\setcounter{equation}{0}

\section{Spaces of coinvariants in terms of filtered tensor product}
\label{sec:3}
In this section we establish a connection between spaces of coinvariants
(\ref{coinv2}) and the filtered tensor product $\fil$ of finite dimensional 
$\g$-modules introduced in \cite{FL}. Motivated by Corollary \ref{CORO},
we restrict our considerations to right ideals $X_1, \cdots, X_N$ 
containing $B_1$, so that $$X_i = S(\gX_i) + B_1$$
where $\gX_i$ is a left ideal of $U(\g)$.
The main statement is Theorem \ref{cnt}.

\subsection{Annihilating ideals}\label{subsec:3.1}
For a $\g$-module $V$, we denote by $\ann V\subset U(\g)$ the kernel
of the structure map $U(\g)\rightarrow \End (V)$.
$\ann V$ is a two-sided ideal and $\ant(\ann V) = \ann V^*$. 
If $V$ is irreducible, then $U(\g)/\ann V\simeq \End (V)\simeq V^*\otimes V$.
The left $\g$-action (resp., the right $\g$-action) on $U(\g)/\ann V$
translates to the left $\g$-action on $V$
(resp., the right $\g$-action on $V^*$).

The representation $\Lkinf$ is a quotient of $U(\gc)$ by the left ideal 
which describes the conditions on the highest weight vector $v^{(k)}_\lambda$.
Hence the space of coinvariants may be viewed as a 
bi-quotient of the algebra $U(\gc)$. 
More precisely, let $e_\theta \in \gn_+$ be the root vector for 
the maximal root $\theta$. 
Introduce the following left ideals of $U(\g[t])$, which annihilate
$v^{(k)}_\lambda$:
\bea
&&\Ila=U(\gc)\otimes_{U(\g)}\gX_\lambda
=U(\gc)\gn_++\sum_{h\in\hh}U(\gc)(h-\br{\lambda,h}),\label{Ila}\\
&&\Ikla=\Ila+U(\gc)(e_{\theta}\otimes t)^{k-\br{\lambda,\hkr}+1}.\label{Ikla}
\ena
These are not enough for the characterization of $v^{(k)}_\lambda$.
However, we have
\begin{lem}\label{LEM5}
We have
\[
U(\g)/(\gX_\lambda+\ann \pi_\lambda)\simeq\pi_\lambda.
\]

In other words, we have the equality of the left ideals (see \eqref{gxl2})
$$\gX'_\lambda=\gX_\lambda+\ann \pi_\lambda.$$
\end{lem}
\begin{proof}
We have
\[
U(\g)/\ann \pi_\lambda\simeq\pi^*_\lambda\otimes\pi_\lambda.
\]
Taking the quotient by the left ideal $\gX_\lambda$ we have
\[
(\pi^*_\lambda\otimes \pi_\lambda)/
(\pi^*_\lambda\gX_\lambda)\otimes \pi_\lambda
\simeq\C v^*_\lambda\otimes\pi_\lambda\simeq \pi_\lambda.
\]
\end{proof}

We have the following characterization of $\Lkinf$ as a quotient of $U(\gc)$.
\begin{prop}\label{prop:3.0}
We have
\be
\Lkinf \cong U(\gc)/(\Ikla+U(\gc)\ann \pi_\lambda).
\en
\end{prop} 
\begin{proof}
Using Lemma \ref{LEM5}, we obtain
\begin{eqnarray*}
\Lkinf&\simeq&(U(\gc)\otimes_{U(\g)}\pi_\lambda)
/U(\gc)(e_{\theta}\otimes t)^{k-\br{\lambda,\hkr}+1}\\
&\simeq&U(\gc)\otimes_{U(\g)}(U(\g)/(\gX_\lambda+\ann \pi_\lambda))
/U(\gc)(e_{\theta}\otimes t)^{k-\br{\lambda,\hkr}+1}\\
&\simeq&U(\gc)/(\Ikla+U(\gc)\ann \pi_\lambda).
\end{eqnarray*}
\end{proof}

Set
\be
\ik = \ann\bigoplus_{\lambda \in \simk} \pi_\lambda.
\en
Clearly $\ik$ is a two-sided ideal and $\ant(\ik) = \ik$.
We have an isomorphism of $U(\g)$-bimodules
\bea
U(\g)/\ik\simeq
\bigoplus_{\lambda\in P_+^{(k)}}\pi_\lambda^*\otimes\pi_\lambda.
\label{iso4}
\ena

For a left ideal $\gX \subset U(\g)$, introduce the 
left $\g$--module
\bea
\pi^{(k)} (\gX) = \bigoplus_{\lambda \in \simk} 
\left(\pi_\lambda^{\gX}\right)^*\otimes \pi_\lambda,
\label{pik}
\ena
and the left ideal
\[
\tilde{\gX}=\gX + \ik\subset U(\g).
\]
\begin{prop}\label{prop:3.1}
We have 
\bea
\pi^{(k)}(\gX) \cong U(\g)/\tilde{\gX}.
\label{iso2}
\ena
In particular, $\pi^{(k)}(\gX)$ is a cyclic $\g$-module, 
the cyclic vector being the sum of canonical vectors
of $\left(\pi_\lambda^{\gX}\right)^*\otimes \pi_\lambda
\cong{\rm Hom}_{\C}(\pi_\lambda^{\gX},\pi_\lambda)$ corresponding to 
the inclusion $\pi_\lambda^{\gX}\subset \pi_\lambda$. 
\end{prop}

\begin{proof}
The isomorphism \eqref{iso2} is clear from \eqref{iso4}:
\[
U(\g)/\tilde{\gX}\cong (U(\g)/\ik)/\gX\cong\bigoplus_{\lambda\in P_+^{(k)}}
(\pi_\lambda^*/\pi_\lambda^*\gX)\otimes\pi_\lambda
\cong\pi^{(k)}(\gX) .
\]
Since $1\in U(\g)$ is mapped in \eqref{iso4} to the canonical vector 
in the right hand side, the statement about the cyclic vector follows.
\end{proof}

The statement of Lemma \ref{LEM5} can be improved as follows.
\begin{cor}\label{COR1}
For $\lambda\in P^{(k)}_+$ we have
\[
U(\g)/(\gX_\lambda+\ik)\cong \pi^{(k)}(\gX_\lambda)\cong \pi_\lambda.
\]
In other words, we have the equality of the left 
ideals
\begin{equation}\label{DEFID}
\gX'_\lambda=\gX_\lambda+\ann \pi_\lambda=
\gX_\lambda+\ik= \widetilde{\gX_\lambda}
\end{equation}
\end{cor}

\subsection{Isomorphism with filtered tensor product}\label{subsec:3.2}

\begin{prop}\label{prop:3.2}
Let $\tilde{\gX}_i \subset U(\g)$ 
be left ideals such that 
$U(\g)/\tilde{\gX}_i$ is finite dimensional, 
$i=1,\cdots,N$. 
Define right ideals of $U(\g[t])$ by 
$\tilde{X}_i=\ant(\tilde{\gX}_i)+B_1$.
Then we have a canonical isomorphism 
\bea
U(\g[t])/\tilde{X}_1\fusn\cdots\fusn \tilde{X}_N(\zz)
\cong
\mathcal{F}_{\zz}\bigl(U(\g)/\tilde{\gX}_1,
\cdots,U(\g)/\tilde{\gX}_N\bigr).
\label{Ufus}
\ena
In this isomorphism the right action of an element $x\in U(\g[t])$
in the left hand side is translated to the left action of $S(x)$
in the right hand side.
\end{prop}

\begin{proof}
Set $U=U(\g[t])$, and let $p_i:U\rightarrow U/\phi_{z_i}(\tilde{X}_i)$ 
be the projection. Consider the composition of maps
\be
U&&
\overset{\Delta^{(N-1)}}{\longrightarrow}
U\otimes\cdots\otimes U
\\
&&\overset{p_1\otimes\cdots\otimes p_N}{\longrightarrow}
U/\phi_{z_1}(\tilde{X}_1)
\otimes\cdots\otimes U/\phi_{z_N}(\tilde{X}_N)
\\&&\overset{\ant\otimes\cdots\otimes \ant}{\longrightarrow}
U/\phi_{z_1}(\tilde{\gX}_1+B_1)\otimes\cdots\otimes 
U/\phi_{z_N}(\tilde{\gX}_N+B_1).
\en
Since $U/\phi_{z_i}(\tilde{\gX}_i+B_1)$ is isomorphic to
the evaluation module $U(\g)/\tilde{\gX}_i$ at $t=z_i$,
the last member is nothing but
\be
\mathcal{F}_{\zz}
(U(\g)/\tilde{\gX}_1,
\cdots,U(\g)/\tilde{\gX}_N).
\en
The map
$U\rightarrow \mathcal{F}_{\zz}
(U(\g)/\tilde{\gX}_1,
\cdots,U(\g)/\tilde{\gX}_N)$
is surjective since the filtered tensor product is cyclic.

Proposition is proved if we show that the kernel of 
$(p_1\otimes\cdots\otimes p_N)\circ\Delta^{(N-1)}$
coincides with $\tilde{X}_1\fusn\cdots\fusn \tilde{X}_N(\zz)$. 
This is shown similarly as in the proof of Proposition \ref{loc}
by using Lemma \ref{LEMTWO}.
\end{proof}

We are now in  a position to state the main result of this section. 
\begin{thm}\label{cnt}
Let $\gX_i \subset U(\g)$ $(i=1,\cdots,N)$ be a set of left ideals.
Define right ideals of $U(\gc)$ by $X_i = \ant(\gX_i) + B_1$.
Then we have a canonical isomorphism of filtered vector space
\begin{equation}
\label{SPC}
\Lkinf/ X_1 \fusn \dots \fusn X_N(\zz) \cong
\fil\left(\pi^{(k)}(\gX_1), \dots , \pi^{(k)}(\gX_N)\right)/
\ant\left(\Ikla\right)
\end{equation}
where $\pi^{(k)}(\gX)$ is defined in \eqref{pik} and the right ideal
$S(\Ikla)$ is given by \eqref{Ikla}$:$
\[
S(\Ikla)=\gn_+U(\gc)+\sum_{h\in\hh}(h+\br{\lambda,h})U(\gc)
+(e_{\theta}\otimes t)^{k-\br{\lambda,\hkr}+1}U(\gc).
\]
\end{thm}

\begin{proof}
Set 
\be
\txx_i = \gX_i + \ik. 
\en
We have
\be
\tx_i{\buildrel{\rm def}\over=}\ant(\txx_i) + B_1=X_i + \ik U(\gc). 
\en
Since $\ik$ annihilates the $\g$-module $\Lkinf/B_1$, 
we have $\Lkinf/ X_i \cong \Lkinf/\tx_i$. 
Comparing dimensions using Theorem~\ref{vermul},
we find that the canonical projection 
\bea
\Lkinf / X_1 \fusn \dots \fusn X_N(\zz) \longrightarrow
\Lkinf / \tx_1 \fusn \dots \fusn \tx_N(\zz)
\label{is1}
\ena
is an isomorphism. 

We have
\be
U(\g[t])/\tilde{X}_1\fusn\cdots\fusn\tilde{X}_N(\zz)
&\cong&
\mathcal{F}_{\zz}(U(\g)/\tilde{\gX}_1,
\cdots,
U(\g)/\tilde{\gX}_N)\\
&&\mbox{(by Proposition \ref{prop:3.2})}\\
&\cong&
\mathcal{F}_{\zz}
\bigl(\pi^{(k)}(\gX_1),
\cdots,
\pi^{(k)}(\gX_N)\bigr).\\
&&\mbox{(by Proposition \ref{prop:3.1})}
\en
Therefore, from Proposition \ref{prop:3.0} follows that  
\bea
&&
\Lkinf/\tilde{X}_1 \fusn \dots \fusn\tilde{X}_N(\zz)
\label{is2}\\
&&\cong 
\mathcal{F}_{\zz}
\bigl(\pi^{(k)}(\gX_1), \dots , \pi^{(k)}(\gX_N)\bigr)
/\left(\ant(I^{(k)}_{\lambda})+\ann \pi_{\lambda^*}U(\gc)\right).
\nn
\ena

The filtered tensor product
$V=\mathcal{F}_{\zz}\bigl(\pi^{(k)}(\gX_1),\cdots,\pi^{(k)}(\gX_N)\bigr)$
is a finite dimensional $\g$-module. 
Hence it has the form $V=\oplus_{\mu} V_\mu$ 
where $V_\mu$ is a direct sum of copies of $\pi_\mu$. 
Noting that $S(I_\lambda)
=\gn_+U(\gc)+\sum_{h\in\hh}(h+\br{\lambda,h})U(\gc)$,
we have
$V/\ant(I_{\lambda})\simeq V_{\lambda^*}/\ant(I_\lambda)$, and 
$\mathop{{\rm Ann}}(\pi_{\lambda^*})$ acts as $0$ on it. 
Therefore 
\bea
V/\left(\ant(I^{(k)}_{\lambda})+\mathop{{\rm Ann}}(\pi_{\lambda^*})U(\gc)\right)
\simeq V/\ant(I^{(k)}_{\lambda}).
\label{is3}
\ena
Combining
\eqref{is1},\eqref{is2} and \eqref{is3}, 
we obtain the assertion. 
\end{proof}

For filtered tensor product of irreducible modules, 
we obtain the following rule for the dimension.
\begin{cor}\label{cor:3.2}
Suppose that $\lambda_1,\ldots,\lambda_N\in P^{(k)}_+$. Then, we have
\bea
\dim \mathcal{F}_{\zz}(\pi_{\lambda_1},\cdots,\pi_{\lambda_N})/
S(I^{(k)}_\lambda)
=\bigl([\lambda_1^*]\cdots[\lambda_N^*]:[\lambda]\bigr)_k.
\label{specialdim}
\ena
\end{cor}
\begin{proof}
Set $\gX_i=\gX_{\lambda_i}$.
{}From Corollary \ref{pik} we have $\pi^{(k)}(\gX_i)\cong\pi_{\lambda_i}$.
Therefore, we have
\[
\mathcal{F}_{\zz}(\pi_{\lambda_1},\cdots,\pi_{\lambda_N})/S(I^{(k)}_\lambda)
\cong\mathcal{F}_{\zz}(\pi^{(k)}(\gX_1),\cdots,\pi^{(k)}(\gX_N))/S(I^{(k)}_\lambda)
\]
Using Theorem \ref{cnt}, and applying Theorem \ref{vermul}
and Corollary \ref{CORO}, we have the assertion.
\end{proof}

\noindent{\bf Example 1.}\quad
For $\lambda_i\in P^{(k)}_+$, set $\gX_i=\gX_{\lambda_i}$ or
$\gX_i=\gX_{\lambda_i}'$
(see \eqref{gxl}, \eqref{gxl2}).
Then, from Corollary~\ref{COR1}, we have
\be
\Lkinf/X_1\fusn\cdots\fusn X_N(\zz)
\simeq \mathcal{F}_{\zz}\left(\pi_{\lambda_1},\cdots,\pi_{\lambda_N}\right)
/S(\Ikla).
\en
By Proposition~\ref{intr}, we have an isomorphism to
the space of conformal coinvariants on $\C P^1$:
\[
\Lkinf/X_1\fusn\cdots\fusn X_N(\zz)\cong
L_\lambda\boxtimes L_{\lambda_1}\boxtimes\dots\boxtimes L_{\lambda_N}
/{\g}^{\C P^1\setminus \{\infty,z_1, \dots,z_N\}}.
\]

\noindent{\bf Example 2.}\quad Let $\gX_i=0$ for all $i$. Then 
$\pi^{(k)}(0)=\sum_{\lambda\in P_+^{(k)}}\pi_\lambda^*\otimes\pi_\lambda$, and 
\be
\Lkinf/B_{1,\zz}\simeq 
\mathcal{F}_{\zz}\left(\pi^{(k)}(0),\cdots,\pi^{(k)}(0)\right)/\ant(\Ikla).
\en
The left hand side has dimension
$\left(\bigl(\sum_{\mu\in P_+^{(k)}}\dim\pi_\mu\cdot [\mu]\bigr)^N:
[\lambda]\right)_k$. 
\medskip

\noindent{\bf Example 3.}\quad Let $\gX_i=\gX=U(\g)\gn_+$ for all $i$. 
We choose the sum of highest weight vectors of $\pi_\lambda$
as cyclic vector of the module
$\pi^{(k)}(\gX)\cong\oplus_{\lambda\in P_+^{(k)}}\pi_\lambda$.
Set $X=\gn_+U(\g)+B_1
=(\gn_+\otimes\C[t]\oplus{\mathfrak h}\otimes t\C[t]\oplus
\gn_-\otimes t\C[t])U(\g[t])$. 
{}From Proposition \ref{prop:2.0}, 
the fusion right ideal is generated by the Lie subalgebra
\[
{\mathfrak a}(\zz)=
\gn_+\otimes\C[t]
\oplus{\mathfrak h}\otimes\prod_{i=1}^N(t-z_i)\C[t]
\oplus\gn_-\otimes\prod_{i=1}^N(t-z_i)\C[t].
\]
Then, we have
\be
\Lkinf/{\mathfrak a}(\zz)
\simeq \mathcal{F}_{\zz}\left(\pi^{(k)}(\gX),\cdots,\pi^{(k)}(\gX)\right)
/S(\Ikla).
\en
The left hand side has dimension 
$\left(\bigl(\sum_{\mu\in P_+^{(k)}} [\mu]\bigr)^N:[\lambda]\right)_k$. 
\medskip

\noindent{\bf Example 4.}\quad Let $\g=\slt$ with standard generators $e,f,h$. 
Let $\pi_l$ $(0\le l\le k)$ denote
the $(l+1)$-dimensional irreducible representation, 
$[l]$ the corresponding element of $\mathcal{V}^{(k)}$,
and $l$ the highest weight.
For $0\le m\le k$, set
\be
&&\gY_m=U(\g)e^{m+1}+U(\g)f,
\\
&&\overline{\gY}_m=U(\g)f^{m+1}+U(\g)e,
\en
and $Y_m=S(\gY_m)+B_1$, $\overline{Y}_m=S(\overline{\gY}_m)+B_1$.
Then 
\be
\pi^{(k)}(\gY_m)= \oplus_{l=0}^k \left(\pi_l^{\gY_m}\right)^*\otimes\pi_l
\cong\oplus_{l=0}^m \pi_l,
\en
with the cyclic vector being the sum of lowest weight vectors. 
Likewise 
$\pi^{(k)}(\overline{\gY}_m)\cong\oplus_{l=0}^m \pi_l$
has the sum of highest weight vectors as cyclic vector. 
Let $Y$ be the fusion right ideal of
\be
\overbrace{Y_1,\cdots,Y_1}^{m_1},
\cdots,
\overbrace{Y_k,\cdots,Y_k}^{m_k},~
\overbrace{\overline{Y}_1,\cdots,\overline{Y}_1}^{n_1},
\cdots,
\overbrace{\overline{Y}_k,\cdots,\overline{Y}_k}^{n_k}.
\en
Then the space of coinvariants $L^{(k)}_l(\infty)/Y$ 
has dimension
\be
\left(([0]+[1])^{m_1+n_1}\cdots([0]+\cdots+[k])^{m_k+n_k}:[l]\right)_k.
\en
For $m=k$, one can replace 
$\gY_k$ by $U(\g)f$ without changing the space of coinvariants (\ref{SPC})
since $\pi^{(k)}(\gY_k)=\pi^{(k)}(U(\g)f)$. 
Similarly, one can replace $\overline{\gY}_k$ by $U(\g)e$.
Therefore, instead of $Y_k$ and $\overline Y_k$, one can use
\begin{eqnarray*}
Z=(e\otimes t\C[t]\oplus h\otimes t\C[t]\oplus f\otimes\C[t])U(\g),\\
\overline Z=(e\otimes\C[t]\oplus h\otimes t\C[t]\oplus f\otimes t\C[t])U(\g),
\end{eqnarray*}
respectively.

Consider the special case where $m_i=n_i=0$ for $1\le i\le k-1$ 
and $m_k=M$, $n_k=N$. In this case,  
the fusion right ideal 
\[
\overbrace{Z\fusn\cdots\fusn Z}^M
\fusn\overbrace{\overline Z\fusn\cdots\fusn\overline Z}^N
(z_1,\ldots,z_M,z'_1,\ldots,z'_N)
\]
is generated by the Lie subalgebra
\be
&&{\mathfrak a}^{(M,N)}(\zz,\zz')
\\
&&
=
e\otimes\prod_{i=1}^M(t-z_i)\C[t]
\oplus h\otimes\prod_{i=1}^M(t-z_i)\prod_{i=1}^N(t-z'_i)\C[t]
\oplus f\otimes\prod_{i=1}^N(t-z'_i)\C[t].
\en
The space of coinvariants 
with respect to ${\mathfrak a}^{(M,N)}(\zz,\zz')$,  
in particular in the degeneration 
limit $z_1=\cdots=z_M=z_1'=\cdots=z_N'=0$, 
is the subject of the works \cite{FKLMM2,FKLMM3}.

\subsection{Fusion Kostka polynomials}\label{subsec:3.4}
In this subsection, we first recall the definition of 
$q$-supernomials, Kostka
polynomials and their generalizations from the literature, and then propose
their counterparts in the theory of fusion product.

For two partitions $\lambda,\mu$ such that $|\lambda|=|\mu|$,
the Kostka-Foulkes polynomials $K_{\lambda,\mu}(q)$ are the entries of the
transition matrix from Schur functions $s_\lambda(x)$ to Hall-Littlewood
functions $P_\mu(x;q)$ (see \cite{Mac} p. 239):
\begin{equation}\label{FRENCH}
s_\lambda(x) = \sum_{\mu} K_{\lambda,\mu}(q) P_{\mu}(x,q).
\end{equation}
If $|\lambda|\neq |\mu|$, $K_{\lambda,\mu}(q)=0$.

{}From the theory of Schur functions, it is known that the Kostka
number $K_{\lambda,\mu}(1)$ is the cardinality of the set of
semi-standard tableaux ${\rm Tab}(\lambda,\mu)$ of shape $\lambda$ and
weight $\mu$.  Therefore it was conjectured by \cite{Fo} that there
should exist a charge
statistic $c: {\rm Tab}(\lambda,\mu) \to \Z_{\geq0}$ on semi-standard
tableaux  such that
$$
K_{\lambda,\mu}(q)=\sum_{t\in Tab(\lambda,\mu)} q^{c(t)}.
$$
This charge statistic was indeed found in \cite{LS}.

In the representation theory of $\mathfrak{gl}_n$, the Kostka number
appears in two ways. In this connection, we have the restriction on
the length of $\lambda$: $\ell(\lambda)\leq n$. First, for $\mu$
satisfying $\ell(\mu)\leq n$, it is equal to the
multiplicity of the weight $\mu$ in the $\mathfrak{gl}_n$ representation
$V_\lambda$ with highest weight $\lambda$.

Secondly, for general $\mu$ such that $\ell(\mu)=N$, it is equal to the
multiplicity of the $\mathfrak{gl}_n$-representation $\pi_\lambda$ in
the tensor product of the $N$ symmetric tensor representations
$V_{(\mu_i)}$ of $\mathfrak{gl}_n$:
\begin{equation}\label{TENS}
\pi_{(\mu_1)}\otimes\cdots\otimes \pi_{(\mu_{N})} \simeq
\oplus_{\lambda:\ell(\lambda)\leq n} K_{\lambda,\mu}(1) \pi_{\lambda}.
\end{equation}

{}From the latter viewpoint, the Kostka polynomial can be regarded as the
$q$-multiplicity of the representation $\pi_\lambda$ in the tensor product
of symmetric tensor representations. Namely, it is the Hilbert polynomial in
the variable $q$ of the graded space
\[
{\rm Hom}_{\gln}(\pi_\lambda,\pi_{(\mu_1)}\otimes\cdots\otimes\pi_{(\mu_N)}),
\]
where the grading is given by the charge statistic.

This interpretation first appeared in
\cite{KR}, where the completeness problem for Bethe ansatz solutions
of the XXX model was considered. This is also the place where 
fermionic formulas for the Kostka polynomials first appeared. In that
approach, the fermionic formulas were proven by constructing a (charge
preserving) bijection between rigged configurations appearing in the
completeness problem for the Bethe equations, and Young tableaux. 
The term ``fermionic formulas'' were invented in
the works of Stony Brook group  \cite{KKMM},\cite{KM}, 
in which it was found that the Bethe ansatz solution 
gives, in the conformal limit, 
the characters of conformal field theory 
in fermionic form, where the
$q$-grading is the linearized quasi-particle 
energy function. 
This
gives a physical 
interpretation to the $q$-grading used in \cite{KR}.

Kostka polynomials were generalized \cite{KS,SW1,SW2} to the case
where the representations 
in the tensor product correspond to highest
weight representations for 
$\sln$ with highest weight $m\lambda_i$  
(rectangular highest weight)  
where $\lambda_i$ $(1\leq i\leq n-1)$ are 
the fundamental weights.  

There are two variants of Kostka polynomials: 
the $q$-supernomials and
the level-restricted Kostka polynomials.

The $q$-supernomials are related to the Kostka polynomials by
\[
S_{\lambda,\mu}(q)=\sum_\nu K_{\nu,\lambda}(1)K_{\nu,\mu}(q).
\]
Namely, it is the $q$-multiplicity of the weight $\lambda$ in the tensor
product (\ref{TENS}). Fermionic formulas for the 
$q$-supernomial 
$S_{\lambda,\mu}(q)$ is known \cite{HKKOTY}, and recursion equations
for the generalized version in the sense of rectangular diagram
are studied \cite{SW2}.

The level restriction has an origin in the corner transfer matrix method
of Baxter. It was systematically studied as one-dimensional configuration sums
in connection with the characters of integrable highest weight representations
of affine Lie algebras. In this context, the level parameter appears
as the level of representation. Kashiwara's theory of crystal base gave a
solid basis to this approach. Quite generally, the crystal base of level $k$
integrable highest weight representations of affine Lie algebras are given in
terms of paths, i.e., elements in semi-infinite product $B^{\otimes\infty/2}$
of a finite dimensional affine crystal $B$. The $q$-statistics 
for paths are
related to the eigenvalues of corner transfer matrix in the zero temperature
limit.

The connection between Kostka polynomials and crystal base
was discovered and generalized in 
\cite{NY,SW1,SW2,HKKOTY}.    
In these works were considered 
finite and inhomogeneous paths,  
i.e., elements in a finite tensor 
product $B_1\otimes\cdots\otimes B_N$
of affine crystals $B_i$.
There are actions of Kashiwara operators $\tilde e_i$, $\tilde f_i$
$(i=0,1,\ldots,r)$ on affine crystals corresponding to the Chevalley
generators $e_i,f_i$ of the affine Lie algebras. 
In this language, the generalized
Kostka polynomial is the $q$-statistical sum over paths restricted by
the highest weight conditions.
\[
\tilde e_i(b_1\otimes\cdots\otimes b_N)=0\quad(i=1,\ldots,r).
\]
The $q$-supernomial is nothing but the unrestricted sum, and the
level $k$ restricted Kostka polynomial is the sum with the integrability
conditions
\[
\tilde f_0^{k-\left<\lambda,\theta^{\vee}\right>+1}
(b_1\otimes\cdots\otimes b_N)=0.
\]

For the generalized Kostka polynomials, 
fermionic formulas are known 
in the case of rectangular weights \cite{KSS}. 
To our knowledge, no formulas are available at present 
for an arbitrary sequence of highest weights. 
Fermionic formulas for the level restricted 
(generalized) Kostka polynomials 
were obtained in \cite{SS} (see also \cite{HKKOTY}).

Now let us formulate the fusion product analogs of the
Kostka and related polynomials. 
In what follows we consider only left ideals
$\gX\subset U(\g)$ which are also $P$-graded,  
\be
\gX=\oplus_{\nu\in P}\gX_\nu,
\en
where
$\gX_\nu=\{x\in\gX\mid [h,x]=\br{\nu,h}x~~~\forall h\in\hh\}$.
The space of coinvariants are then $P$-graded,  
as well as filtered by the degree of $t$. 
In the isomorphism \eqref{Ufus} of Proposition \ref{prop:3.2},
the left action by $h\in \mathfrak{h}$ on the left hand side
is translated to the right action by $S(h)=-h$ on 
the filtered tensor product (the right hand side). 
The latter commutes with the action of $\g$,
and is not to be confused with the left action of $h\in \g$. 
To be concrete, let us 
take a basis $\{u_{\lambda,i}\}$ of $\pi_\lambda^{\gX}$ and 
the dual basis $\{u^i_{\lambda}\}$ of $(\pi_\lambda^{\gX})^*$.
The cyclic vector of $\pi^{(k)}(\gX)$ is 
$u^{(k)}(\gX)
=\sum_{\lambda\in P^{(k)}_+}\sum_iu^i_{\lambda}\otimes
u_{\lambda,i}$, 
and the right action of $S(h)$ $(h\in\mathfrak{h})$ 
on a vector $xu^{(k)}(\gX)$ $(x\in\g)$ is given by 
\bea
xu^{(k)}(\gX)S(h)=-
\sum_{\lambda\in P^{(k)}_+}\sum_i(u^i_{\lambda}h)\otimes
(xu_{\lambda,i}).  
\label{Pgr}
\ena
The $P$-grading on the space of coinvariants
\eqref{SPC}
of Theorem \ref{cnt} is induced from the above $P$-grading. 

In general, let $W=\oplus_{\nu\in P}W_\nu=\varinjlim F^j$
be a $P$-graded vector space equipped with a filtration 
$0=F^{-1}\subset F^0\subset F^1\subset\cdots$ 
by $P$-graded subspaces $F^j$.
We define the character of the associated graded space by 
\be
{\rm ch}\gr W=
\sum_{j\in\Z_{\ge 0},\nu\in P}
\dim (F^j/F^{j-1})_\nu\,q^je^\nu,
\en
where $e^\nu$ stands for the formal exponential symbol.
When the $P$-grading is trivial, we use also the notation ${\rm ch}_q$. 

The foregoing discussions  
lead us to consider three kinds of polynomials in $q$,    
$\mathcal{S}_{\lambda}(V_1,\cdots,V_N)$, 
$\fk_\lambda (V_1, \dots, V_N)$ 
and $\fk^{(k)}_\lambda (V_1, \dots, V_N)$. 
\medskip

\noindent{\bf Fusion $q$-supernomial coefficient}\quad
\bea
&&
\ch~\mathop{\rm gr} 
\mathcal{F}_{\mathcal{Z}}(V_1, \cdots, V_N)
=\sum_{\lambda\in P}e^\lambda 
\mathcal{S}_{\lambda}(V_1,\cdots,V_N).
\label{fussup}
\ena
\medskip

\noindent{\bf Fusion Kostka polynomial}\footnote{The fusion Kostka 
polynomial defined in \cite{Kir} is different from ours.}
\quad
\bea
&&
\fk_\lambda (V_1, \dots, V_N) 
=\ch~ \mathop{\rm gr}\Bigl(\mathcal{F}_{\mathcal{Z}}(V_1, \cdots, V_N)/
\ant(\Ila)\Bigr),
\label{fusKos}
\ena
\medskip

\noindent{\bf Fusion restricted Kostka polynomial}\quad
\bea
&&
\fk^{(k)}_\lambda (V_1, \dots, V_N) 
=\ch~\mathop{\rm gr}
\left(\mathcal{F}_{\mathcal{Z}}(V_1, \cdots, V_N)/\ant(\Ikla)\right).
\label{fusresKos}
\ena
\medskip

The first two are simply related by an alternating sum over the Weyl group as 
\be
\fk_\lambda (V_1, \dots, V_N) 
= 
\sum_{w \in W}
(-1)^{l(w)} 
\mathcal{S}_{w\circ \lambda} (V_1, \dots, V_N).
\en
It is also clear that $\fk_\lambda (V_1, \dots, V_N)$ is obtained from 
$\fk^{(k)}_\lambda(V_1, \dots, V_N)$ in the limit $k \to \infty$.

We expect that, in the case when 
$V_1,\cdots,V_N$ are tensor powers of the fundamental representations
discussed in 
\cite{HKKOTY},\cite{SW2}, 
the `fusion-product' version defined above reduces 
to the corresponding object known in the literature.
We expect also the following relation:
\medskip

\noindent{\bf Conjecture.}\quad
\bea
\phantom{conj}
\fk^{(k)}_\lambda (V_1, \dots, V_N) = 
\sum_{{w \in W_{\rm aff}}\atop{w\circ\lambda \in P_+}} 
(-1)^{l(w)} q^{d_w(\lambda,k)}
\fk_{w\circ \lambda} (V_1, \dots, V_N).
\label{mainconj}
\ena
The quantity $d_w(\lambda,k)$
is defined in (\ref{shift}).
For the cases corresponding to Examples 2,3, 
\eqref{mainconj} was  put forward in \cite{FL} 
as the `main conjecture'.

\setcounter{section}{3}
\setcounter{equation}{0}

\section{The case $\g=\slt$}\label{sec:4}
In this section we restrict ourselves to the 
case $\g=\slt$ as mentioned in \eqref{sl2case}. 
In this special case, 
we shall identify the fusion-product version of 
Kostka/level-restricted Kostka polynomials with the 
ordinary ones and prove the identity \eqref{mainconj}. 
As before, we write $\pi_l$ for the $(l+1)$-dimensional irreducible 
representation and $[l]$ for the corresponding element of $\mathcal{V}^{(k)}$. 
We also write the left ideal \eqref{Ikla} as $I^{(k)}_l$, 
and $\ant(I^{(k)}_l)$ as $\br{h[0]+l,e[0],e[1]^{k-l+1}}$, 
where $x[i]=x\otimes t^i$.  

\subsection{level-restricted Kostka polynomials}\label{subsec:4.1}
First we recall a few basic facts 
about Kostka and level-restricted Kostka polynomials for $\slt$, 
following the exposition of \cite{SS}. 
We fix the notation as follows. 
Let $\mathbf{m}=(m_1,\cdots,m_k)$ 
be a $k$-tuple of non-negative integers, and let $0\le l\le k$.  
Here and after, we use boldface letters to represent a vector. 
The transposition symbol is omitted. 
Set 
\be
&&K_{l,\mathbf{m}}(q)=
q^{\Vert \mathbf{m}\Vert}K_{\lambda\, R(\mathbf{m})}(q^{-1}),
\\
&&K_{l,\mathbf{m}}^{(k)}(q)=
q^{\Vert \mathbf{m}\Vert}K^k_{\lambda\, R(\mathbf{m})}(q^{-1}).
\en
Here the right hand sides stand for the Kostka and level-restricted
Kostka polynomials in the notation of \cite{SS}, 
\be
&&\lambda=\left(\frac{|\mathbf{m}|+l}{2},
\frac{|\mathbf{m}|-l}{2}\right),\\
&&R(\mathbf{m})=(k^{m_k},\cdots,2^{m_2},1^{m_1}), 
\en
and
$|\mathbf{m}|=\sum_{i=1}^kim_i$, 
$2\Vert \mathbf{m}\Vert=\sum_{1\le i,j\le k}\min(i,j)m_im_j-|\mathbf{m}|$.
Set further
\bea
&&A_{ab}=\min(a,b), 
\quad
v_a=\max(a-k+l,0).
\label{A}
\ena
We use the $q$-binomial symbol
\be
&&
\qbin{\mathbf{m}}{\mathbf{n}}
=\prod_{i=1}^k\qbin{m_i}{n_i},
\\
&&\qbin{m}{n}=\begin{cases}
\displaystyle\frac{[m]!}{[n]![m-n]!} & (0\le n\le m),\\
0 & \mbox{otherwise},\\
\end{cases}
\en
where $[n]!=\prod_{j=1}^n\bigl((1-q^j)/(1-q)\bigr)$. 
We shall use the following facts.
\begin{description}
\item[Fermionic formula] 
\bea
K_{l,\mathbf{m}}^{(k)}(q)
=
\sum_{\mathbf{s}\in\Z_{\ge 0}^k
\atop 2|\mathbf{s}|=|\mathbf{m}|-l}
q^{\mathbf{s}A\mathbf{s}+\mathbf{v}\mathbf{s}}
\qbin{A(\mathbf{m}-2\mathbf{s})-\mathbf{v}+\mathbf{s}}
{\mathbf{s}}.
\label{RKos}
\ena
\item[Alternating sum formula]
\bea
K_{l,\mathbf{m}}^{(k)}(q) 
&=&\sum_{i\ge 0}q^{(k+2)i^2+(l+1)i}
K_{2(k+2)i+l,\mathbf{m}}(q) 
\label{bos}\\
&&-\sum_{i> 0}q^{(k+2)i^2-(l+1)i}
K_{2(k+2)i-l-2,\mathbf{m}}(q).  
\nonumber
\ena
\item[Verlinde number]
In the Verlinde algebra $\mathcal{V}^{(k)}$, we have 
\bea
[1]^{m_1}\cdots [k]^{m_k}
=\sum_{l=0}^k K_{l,\mathbf{m}}^{(k)}(1) [l].
\label{Number}
\ena
\end{description}
Eq. \eqref{RKos} and \eqref{bos} are taken from Theorem 6.2 and eq. (6.8) 
in \cite{SS}, respectively (after simplifications due to $n=2$).
In \cite{HKKOTY}, $K_{l,\mathbf{m}}^{(k)}(q)$ is 
defined as a generating sum over 
level $k$ restricted paths in affine crystals. 
With this definition, \eqref{Number} follows from the fact that 
the structure constants of the Verlinde algebra 
are the same as the number of successive triples 
in a path satisfying the level-restriction conditions. 

\subsection{The goal}\label{subsec:4.2}
The fusion product is defined to be 
the associated graded vector space of the filtered tensor product 
\be
V_1*\cdots*V_N(\zz)=\gr\mathcal{F}_{\zz}(V_1,\cdots,V_N).
\en
It gives a meaning to the 
`limit' of the filtered
tensor product where all points $z_i$ tend to $0$, 
as explained in Appendix. 
We consider the filtered tensor product and fusion product
of cyclic modules $V_i=\pi_{l_i}$, 
choosing the highest weight vector as cyclic vector.
Let $\mathbf{m}=(m_1,\cdots,m_k)$ and set 
\bea
V_{\mathbf{m}}=
\overbrace{\pi_k*\cdots*\pi_k}^{m_k}
*\cdots*
\overbrace{\pi_1*\cdots*\pi_1}^{m_1}
(\zz).
\label{fus2}
\ena
It is known (for $\mathfrak{g}=\mathfrak{sl}_2$) that 
\eqref{fus2} is independent of $\zz$ \cite{FL,FF}. 
Our goal is to show the following.
\begin{thm}\label{thm:4.1}
(i)Notations being as above, we have
\be
\ch_q V_{\mathbf{m}}/\br{h[0]+l,e[0],e[1]^{k-l+1}}
=K_{l,\mathbf{m}}^{(k)}(q).
\en
(ii)The identity \eqref{mainconj} is true in this case. 
\end{thm}

In the next subsection, we prove the inequality
\bea
\ch_q V_{\mathbf{m}}/\br{h[0]+l,e[0],e[1]^{k-l+1}}
\le K_{l,\mathbf{m}}^{(k)}(q).
\label{estimate}
\ena
Namely, each coefficient of powers in $q$ in the right hand side is
greater than or equal to that in the left hand side.

Assuming \eqref{estimate}, let us prove Theorem \ref{thm:4.1}.
\medskip

{\it Proof of Theorem \ref{thm:4.1}}.\quad 
{}From Corollary \ref{cor:3.2} and 
\eqref{jump}, we have for any $V_i=\pi_{l_i}$ that
\be
([l_1]\cdots[l_N]:[l])_k
&=&
\dim \mathcal{F}_{\zz}(V_1,\cdots,V_N)/\br{h[0]+l,e[0],e[1]^{k-l+1}}
\\
&&
\le \dim V_1*\cdots*V_N/\br{h[0]+l,e[0],e[1]^{k-l+1}}.
\en
Combining this with \eqref{Number}, we obtain 
\be
K_{l,\mathbf{m}}^{(k)}(1)\le 
\dim V_{\mathbf{m}}/\br{h[0]+l,e[0],e[1]^{k-l+1}}.
\en
Therefore, it follows from \eqref{estimate} that equality takes place
in all intermediate steps. In particular, assertion (i) is implied by 
the equality in \eqref{estimate}. Assertion (ii) is then 
an immediate consequence of the known alternating sum formula \eqref{bos}
for level-restricted Kostka polynomials.
\qed

\subsection{Gordon filtration}\label{subsec:4.3}
It remains to prove the inequality \eqref{estimate}.
{}For this purpose we invoke the `functional' description 
of the dual space \cite{FS}.

It was shown in \cite{FF} that, 
as graded vector space, the dual of $V_{\mathbf{m}}$
is isomorphic to the space of all symmetric polynomials 
$f(x_1,\cdots,x_s)$ $(s\ge 0)$ satisfying the condition
\bea
&&\deg_x f(\overbrace{x,\cdots,x}^{a},x_{a+1},\cdots,x_s)
\le M_a-a \label{degf}\\
&&\qquad\qquad 
\qquad \mbox{for any }~ a=1\cdots,s, 
\nonumber
\ena
where
\be
M_a=\sum_{i=1}^k\min(a,i)m_i.
\en
Therefore the dual to the quotient space of \eqref{fus2} by
$\br{h[0]+l,e[0],e[1]^{k-l+1}}$ is realized as follows. 

\begin{lem}\label{lem:4.2}
As a graded vector space, 
$(V_{\mathbf{m}}/\br{h[0]+l,e[0],e[1]^{k-l+1}})^*$ 
is isomorphic to the space $\F$ 
of all symmetric polynomials satisfying \eqref{degf} as well as the conditions
\begin{enumerate}
\item The number of variables is 
$s=\frac{1}{2}(|\mathbf{m}|-l)$, 
\item $f(0,x_2,\cdots,x_s)=0$,
\item If $s\ge k-l+1$, then 
setting 
$f(x_1,\cdots,x_s)
=\bigl(\prod_{i=1}^s x_i\bigr)g(x_1,\cdots,x_s)$ we have
\bea
&&g(\overbrace{0,\cdots,0}^{k-l+1},x_{k-l+2},\cdots,x_s)=0.
\label{fg}
\ena
\end{enumerate}
\end{lem}

\begin{lem}\label{lem:4.1}
Assume that $s\ge k+1$. 
{}For any $f\in\F$ we have 
\be
f(\overbrace{x,\cdots,x}^{a},x_{a+1},\cdots,x_s)=0
\quad ({\rm if}~~a\ge k+1). 
\en
\end{lem}
\begin{proof}
Let $g(x_1,\cdots,x_s)$ be as in Lemma \ref{lem:4.2} (iii).
{}From \eqref{fg}, we find 
\bea
&&\bigl(\frac{\partial}{\partial x}\bigr)^j
g(\overbrace{x,\cdots,x}^{a},x_{a+1},\cdots,x_s)\bigl|_{x=0}=0 \\
&&
\qquad \qquad \qquad (0\le j\le l,~~a\ge k+1). 
\nonumber
\ena
Namely we have 
\begin{equation}\label{g}
g(\overbrace{x,\cdots,x}^{a},x_{a+1},\cdots,x_s)=O(x^{l+1}). 
\end{equation}
We prove the assertion by descending induction on $a$. 
{}First let $a=s$. Since $s\ge k+1$, we have 
$M_a=|\mathbf{m}|$. Then \eqref{degf} and Lemma \ref{lem:4.2} (i) imply
$\deg_x g(x,\cdots,x)\le |\mathbf{m}|-2s=l$. 
Hence we have $g(x,\cdots,x)=0$ by \eqref{g}. 
Suppose the assertion is true for $a+1$ $(a\ge k+1)$. 
Then $g(x_1,\cdots,x_s)$ can be written as 
\be
g(\overbrace{x,\cdots,x}^{a},x_{a+1},\cdots,x_s)
=\prod_{i=a+1}^s(x-x_i)^2 h(x;x_{a+1},\cdots,x_s).
\en
{}From \eqref{degf} we find $\deg_x h\le l$. Therefore 
\eqref{g} implies $h=0$, 
and the proof is over. 
\end{proof}

In order to estimate the character of $\F$, 
we apply the standard technique of Gordon filtration 
\cite{FS}. 
Let $\lambda$ be a partition of width $\lambda_1\le k$, and write 
$\lambda=(k^{s_k}\cdots2^{s_2}1^{s_1})$.
{}For a symmetric polynomial $f$ we define 
\be
&&\varphi_\lambda(f)(x^{(1)}_1,\cdots,x^{(1)}_{s_1};\cdots;
x^{(k)}_1,\cdots,x^{(k)}_{s_k})
\\
&&=f(x^{(1)}_1,\cdots,x^{(1)}_{s_1},
\cdots,
\overbrace{x^{(k)}_1,\cdots,x^{(k)}_1}^{k},\cdots,
\overbrace{x^{(k)}_{s_k},\cdots,x^{(k)}_{s_k}}^{k}),
\en
where each variable $x^{(a)}_i$ appears $a$ times in the right hand side. 
Introduce the lexicographic ordering $\mu>\lambda$ by 
$\mu_1=\lambda_1,\cdots,\mu_{i-1}=\lambda_{i-1}$ and 
$\mu_i>\lambda_i$ for some $i$,  
and set 
\be
\F_{\lambda}
=\bigcap_{\mu>\lambda}\mathop{\rm Ker}\varphi_\mu.
\en
Notice that $\F_{(k)}=\F$ by Lemma \ref{lem:4.1}.
The family of subspaces $\{\F_\lambda\}_{\lambda}$ 
defines a filtration of $\F$, 
whose associated graded space has 
the $\lambda$-th component 
isomorphic to $\varphi_\lambda(\F_\lambda)$. 
An element of the latter has the form 
\be
&&\bigl(\prod_{1\le i\le s_a\atop 1\le a\le k}x^{(a)}_i\bigr)^{a+(a-k+l)_+}
\prod_{1\le i<j\le s_a\atop 1\le a\le k}(x^{(a)}_i-x^{(a)}_j)^{2a}
\\
&&\times 
\prod_{1\le i\le s_a,1\le j\le s_b\atop 1\le a<b\le k}(x^{(a)}_i-x^{(b)}_j)^{2\min(a,b)}
\cdot h,
\en
with a partially symmetric polynomial
\be
h\in 
\C[x^{(1)}_1,\cdots,x^{(1)}_{s_1},\cdots,x^{(k)}_1,\cdots,x^{(k)}_{s_k}]
^{\Sym_{s_1}\times\cdots\times \Sym_{s_k}},
\en
$\Sym_s$ denoting the symmetric group on $s$ letters.
We have the relations for the degree
\be
&&\deg_{x^{(a)}_i}h\le 
\left(A(\mathbf{m}-2\mathbf{s})-\mathbf{v}\right)_a,
\\
&&\mathop{\rm tot.}\deg f=
\mathop{\rm tot.}\deg h+\mathbf{s}A\mathbf{s}+\mathbf{v}\mathbf{s},
\en
where we used the vector notation \eqref{A}, and 
$\mathop{\rm tot.}\deg$ signifies the total degree.
We thus obtain for the character $\ch_q\gr \F$
an upper bound 
\be
\sum_{\mathbf{s}\in\Z_{\ge 0}^k
\atop 2|\mathbf{s}|=|\mathbf{m}|-l}
q^{\mathbf{s}A\mathbf{s}+\mathbf{v}\mathbf{s}}
\qbin{A(\mathbf{m}-2\mathbf{s})-\mathbf{v}+\mathbf{s}}
{\mathbf{s}},
\en
which coincides with the fermionic form
of the level-restricted Kostka polynomial \eqref{RKos}.   
This completes the proof of \eqref{estimate}.

\bigskip

\appendix
\section{Degeneration of coinvariants and fusion product}

Here we give a motivation for the notion of fusion product introduced in
\cite{FL} and used in Section~4.

In this paper we have considered 
fusion of congruence right ideals of $U(\g[t])$ attached to distinct points.
They are parameterized by the configuration space of complex points. 
In general, it is not clear what is the analog of fusion when points collide,
but we can extend the definition to the boundary of the
configuration space ${\bf C}^N\backslash\bigcup\Delta_{[a,b]}$
(see \cite{FM}) in the following way.

For $a\in\C\setminus\{0\}$ we define an automorphism $l_a:\g[t]\to\g[t]$
multiplying the variable $t$ by $a^{-1}$. Clearly we have
$$l_a (X_1 \fusn \dots \fusn X_N(\zz)) = X_1 \fusn \dots \fusn 
X_N(a\zz).$$
These spaces form a family over $\C\setminus \{ 0 \}$, and we want to
extend it to the origin.

\begin{defn}\label{TOP}
{For a subspace $W'$ of a graded vector space
$W=\oplus_{j\ge 0} W_j$,
we denote by ${W'}^{\rm top}$ the graded subspace of $W$ spanned
by homogeneous elements $w\in W_r$,
such that for some $w_j\in W_j$ $(0\le j\le r)$ we have
$w=w_r$ and $w_r+w_{r-1}+\cdots+w_0\in W'$. }
\end{defn}

\def \top {{\rm top}}

Note that $\g[t]$ and therefore $U(\g[t])$ has a natural grading by degree
of $t$. 
Hence, for any subspace $X \subset U(\g[t])$ 
we can consider the subspace $X^\top \subset U(\g[t])$. 
In particular, the space $\left(X_1 \fusn \dots\fusn X_N(\zz)\right)^\top$
can be thought of as the limit of
$X_1 \fusn \dots \fusn X_N(a\zz)$
when the parameter $a$ tends to $0$. 

More generally, combining 
the operation of $\fusn$ and taking $(\dots)^\top$ in 
various ways, 
one can extend the definition of fusion to the compactification
$\C[N]$ of the configuration space in the sense of \cite{FM}.
For example, corresponding to Figure 4 in \cite{FM}, we have
the following fusion of right ideals:
\begin{eqnarray*}
&&Y_1=(X_1\fusn X_3\fusn X_8(z_1,z_3,z_8))^\top,\\
&&Y_2=(X_4\fusn X_6(z_4,z_6))^\top,\\
&&Y_3=(X_2\fusn X_5(z_2,z_5))^\top,\\
&&Y_4=(Y_1\fusn Y_2(z_{138},z_{46}))^\top,\\
&&Y_5=(Y_4\fusn Y_3\fusn X_7(z_{13468},z_{25},z_7))^\top.
\end{eqnarray*}

\medskip

{\bf Question 1.}\  When is the ideal $\left(X_1 \fusn \dots
\fusn X_N(\zz)\right)^\top$ independent of $\zz$? 
More generally, when do the extended fusion right ideals
depend only on the image of the point in
$\C[N]$ under the natural map $\C[N]\to\C^N$?

In the following situation the top part of the right ideal
is given rather explicitly.
\begin{prop}
\label{CITE}
Let ${\mathfrak b}\subset\g[t]$ be a Lie subalgebra. Then, we have
\begin{equation}\label{EQRIGHT}
{\mathfrak b}^{\rm top}U(\g[t])=({\mathfrak b}U(\g[t]))^{\rm top}.
\end{equation}
\end{prop}
\begin{proof}
Let $\C[t]_{\leq d}$ be the space of polynomials of degree less than or
equal
to $d$. Set ${\mathfrak b}_{\leq d}={\mathfrak b}\cap\C[t]_{\leq d}$.
One can construct a basis of ${\mathfrak b}$, $\{B_i\}$,
in such a way that for each $d$ the subset $\{B_i\}\cap{\mathfrak b}_{\leq d}$
constitutes a basis of ${\mathfrak b}_{\leq d}$. Then, $\{B_i^{\rm top}\}$
is a basis of ${\mathfrak b}^{\rm top}$. One can construct a homogeneous
basis
of $\g[t]$ consisting of $\{B_i^{\rm top}\}$ and new members $\{A_i\}$.
Since ${\mathfrak b}$ is a Lie
subalgebra, by Lemma \ref{LEMONE}
$B_{i_1}\cdots B_{i_m}A_{j_1}\cdots A_{j_n}$
where $m\geq1,i_1\leq\cdots\leq i_m$ and $n\geq0,j_1\leq\cdots\leq j_n$
constitute a basis of the right ideal ${\mathfrak b}U(\g[t])$.
Note that
\[
(B_{i_1}\cdots B_{i_m}A_{j_1}\cdots A_{j_n})^{\rm top}
=B_{i_1}^{\rm top}\cdots B_{i_m}^{\rm top}A_{j_1}\cdots A_{j_n}
\]
constitute a basis of $({\mathfrak b}U(\g[t]))^{\rm top}$.
Since ${\mathfrak b}^{\rm top}$ is also a Lie subalgebra,
these elements constitute a basis of ${\mathfrak b}^{\rm top}U(\g[t])$.
Therefore, we have the equality of the right ideals (\ref{EQRIGHT}).   
\end{proof}

In particular, we have a 
positive answer to Question~1 for $X_i=B_{n_i}$, 
since
$$\left(B_{n_1} \fusn \dots \fusn B_{n_N}\right)^{\rm top} 
= B_{n_1 + \dots +n_N}.
$$

Now let us proceed to coinvariants. 
By the definition, there is a canonical 
isomorphism of graded vector space
\begin{equation}\label{GVS}
\gr(W/W')\simeq W/{W'}^{\rm top}.
\end{equation}
On the other hand, if $X$ is a right ideal of a graded algebra $U$ acting
on $W$, the dimension of the space of coinvariants may jump in general.
\bea
{\rm dim}\,W/XW\leq{\rm dim}\,W/X^{\rm top}W.
\label{jump}
\ena

\begin{prop}
There is a surjective map
\begin{equation}
\label{jumpmap}
{L^{(k)}_\lambda(\infty)/\left(X_1\fusn\dots\fusn X_N(\zz)\right)^\top\to 
\gr L^{(k)}_\lambda(\infty)/\left(X_1\fusn\dots\fusn X_N(\zz)\right)}
\end{equation}
\end{prop}

\medskip

{\bf Question 2.}\ When is 
the map \eqref{jumpmap} an isomorphism? 
More generally, 
when is the dimension of coinvariants independent of the point in $\C[N]$?

The answer to this question for $X_i = B_{n_i}$ is positive 
for $\g =sl_2$ (see \cite{FKLMM1}).  
On the other hand, if $\g = E_8$, it is negative even for $X_i=B_1$ with 
$k=1$. This is because for level $1$ the only integrable representation
is the vacuum representation, i.e., the highest weight is $0$. Therefore,
the right hand side is always one-dimensional, while the left hand side
becomes infinite in the limit $N\rightarrow\infty$.

Similar considerations for the filtered tensor product of modules 
can be found in \cite{FL}.
Recall that the associated graded vector space of the filtered tensor
product
\begin{equation}\label{FP}
V_1*\cdots*V_N(\zz)=\gr\mathcal{F}_{\zz}(V_1,\cdots,V_N)
\end{equation}
is called the fusion product. It gives a meaning to the
similar limit of the filtered
tensor product where all points $z_i$ tend to $0$. In \cite{FL} it was
conjectured that, for finite-dimensional cyclic modules over a simple
Lie algebra $\g$, the fusion product $V_1*\cdots*V_N(\zz)$
is independent of the choice of $\zz$.
This conjecture was proved in \cite{FL} in the simplest case where
\bea
&&\g=\slt,\quad V_i:\hbox{ irreducible module}\label{sl2case}\\
&&\mbox{with highest weight vector as cyclic vector},
\nn
\ena
and in \cite{FF} by a different method. Indeed, the fusion product
can also be extended to $\C[N]$ and it can be proved by the
method of \cite{FF} that in this case it depends only on the image in
$\C^N$.

An advantage of this definition is that we can split the map
\eqref{jumpmap} into two.

\begin{thm}\label{Thm:jmp}
\begin{enumerate}
\item The map \eqref{jumpmap} is a composition of the following surjective 
maps
\begin{equation}\label{jumpmap1}
L^{(k)}_\lambda(\infty)/(X_1\fusn\cdots\fusn X_N(\zz))^{\rm top}
\to
\pi^{(k)}(\gX_1)\ast\cdots\ast\pi^{(k)}(\gX_N)/S(\Ikla)
\end{equation}

\begin{equation}\label{jumpmap2}
\pi^{(k)}(\gX_1)\ast\cdots\ast\pi^{(k)}(\gX_N)/S(\Ikla) \to \gr
L^{(k)}_\lambda(\infty)/(X_1\fusn\cdots\fusn X_N(\zz))
\end{equation}

\item Assume that $\gX_i\supset \ik$. Then the map \eqref{jumpmap1} is an
isomorphism.
\end{enumerate}
\end{thm}

\begin{proof}
To prove (i) we just follow the proof of Theorem \ref{cnt}.

\def \tX {\tilde{\gX}}
\def \tx {\tilde{X}}

Let $\tX_i = \gX_i + \ik$. As usual, set $X_i = B_1 + \ant(\gX_i)$, $\tx_i
= B_1 + \ant(\tX_i)$.
Then the natural projection
$$L^{(k)}_\lambda(\infty)/(X_1\fusn\cdots\fusn X_N(\zz)) \to
L^{(k)}_\lambda(\infty)/(\tx_1\fusn\cdots\fusn \tx_N(\zz))$$
is an isomorphism.

So we can split the map~\eqref{jumpmap} into composition of
the projection
$$L^{(k)}_\lambda(\infty)/(X_1\fusn\cdots\fusn X_N(\zz))^{\rm top}
\to L^{(k)}_\lambda(\infty)/(\tx_1\fusn\cdots\fusn \tx_N(\zz))^{\rm
top}$$
and the map~\eqref{jumpmap} for $\tx_i$:
$$L^{(k)}_\lambda(\infty)/(\tx_1\fusn\cdots\fusn \tx_N(\zz))^{\rm top} \to
\gr L^{(k)}_\lambda(\infty)/(\tx_1\fusn\cdots\fusn \tx_N(\zz)).$$

To prove (i)  it remains to
identify $L^{(k)}_\lambda(\infty)/(\tx_1\fusn\cdots\fusn \tx_N(\zz))^{\rm top}$
with the space $\pi^{(k)}(\gX_1)\ast\cdots\ast\pi^{(k)}(\gX_N)/S(\Ikla)$.

We have
\begin{eqnarray*}
&&L^{(k)}_\lambda(\infty)/(\tx_1\fusn\cdots\fusn \tx_N(\zz))^{\rm top}\\  
\simeq&& U(\g[t])/(\Ikla+U(\g[t]){\rm Ann}\pi_\lambda)
/(\tx_1\fusn\cdots\fusn \tx_N(\zz))^{\rm top}\\
&&\hbox{\qquad(by Proposition \ref{prop:3.0})}\\
\simeq &&U(\g[t])
/(\tx_1\fusn\cdots\fusn \tx_N(\zz))^{\rm top}
/(\Ikla+U(\g[t]){\rm Ann}\pi_\lambda)\\
&&\hbox{\qquad(as a bi-quotient by the left and right ideals)}\\
\simeq&& {\rm gr}\Bigl( U(\g[t])
/\tx_1\fusn\cdots\fusn \tx_N(\zz)\Bigr)
/(\Ikla+U(\g[t]){\rm Ann}\pi_\lambda)\\
&&(\hbox{by (\ref{GVS}) and using that }\Ikla\hbox{ is homogeneous})\\
\simeq&& \pi^{(k)}(\gX_1)\ast\cdots\ast\pi^{(k)}(\gX_N)
/(S(\Ikla)+{\rm Ann}\pi_{\lambda^*}U(\g[t]))\\  
&&\hbox{(by Proposition \ref{prop:3.2} 
and by (\ref{FP}))}\\
\simeq &&\pi^{(k)}(\gX_1)\ast\cdots\ast\pi^{(k)}(\gX_N)
/S(\Ikla).\\
&&\hbox{(by repeating the last part of the proof of Theorem \ref{cnt})}
\end{eqnarray*}
To prove (ii) just note that $\gX_i\supset \ik$ implies $X_i = \tx_i$, so
the first map is an isomorphism.
\end{proof}

Let us illustrate how it works for the case of conformal blocks. 
For simplicity suppose $\g = \slt$. 

For this purpose we need ideals $\gX \subset U(\g)$ such that
$\pi^{(k)}(\gX) = \pi_l$. Actually, we can choose them in different
ways. Following \eqref{gxl} and \eqref{gxl2} set
\be
\gX_l=  
U(\slt) e + U(\slt)(h-l),
\en
\be
\gX_l'=
U(\slt) e + U(\slt)(h-l)+ U(\slt)f^{l+1}.
\en

Then $\gX_l' \supset \ik$ because $\gX_l'$ is the annihilating
ideal of the highest weight vector in $\pi_l$ where $0\leq l\leq k$.
Therefore, by
Theorem~\ref{Thm:jmp} the map~\eqref{jumpmap1} is an isomorphism. In
Section~4 we prove that the map~\eqref{jumpmap2} is an isomorphism, so the
composition map~\eqref{jumpmap} is an isomorphism, i.e., the answer to the
Question~2 is positive in this case.

But this is not so for $\gX_l$. Namely, let us show that the
map~\eqref{jumpmap1} can not always be an isomorphism.

Let $X_i = \ant(\gX_{l_i}) + B_1$. In a way similar to the proof of 
Proposition \ref{CITE},
one can show that
$X_1\fusn\cdots\fusn X_N(\mathcal{Z})$
is generated by $e\otimes t^j$ and
$h\otimes t^j-\sum_{i=1}^Nl_iz_i^j$ ($j\ge 0$).
Hence $(X_1\fusn\cdots\fusn X_N(\mathcal{Z}))^{\rm top}$
is generated by
$e\otimes t^j$ ($j\ge 0$),
$h\otimes t^j$ ($j> 0$) and
$h-\sum_{i=1}^Nl_i$.
This means that the left hand side of \eqref{jumpmap1} depends only on
$\sum_{i=1}^Nl_i$, 
but it is easy to see (e.g. by taking $k \to \infty$) that the dimension
of the right hand side depends on each $l_i$.
\bigskip

\noindent
{\it Acknowledgments.}\quad 
This work is partially supported by
the Grant-in-Aid for Scientific Research 
(A1) no.13304010 and (B2) no.14340040, Japan Society for 
the Promotion of Science. BF is partially supported by 
grants RFBR 02-01-01015 and
INTAS-00-00055. 
The work of SL is partially 
supported by the grant RFBR-01-01-00546.
The last stage of this work was carried out 
while the authors were visiting Mathematical Sciences 
Research Institute, Berkeley, March 2002. 


\begin{thebibliography}{[FKLMM3]}

\bibitem[BKMW]{Kir}
L.~Begin, A.N.~Kirillov, P.~Mathieu and M.A.~ Walton,
\newblock
Berenstein-Zelevinski triangles, elementary coupling and fusion rules,
\newblock
{\em Lett. Math. Phys.}, {\bf 28} (1993), 257-268.


\bibitem[FF]{FF}
B.~L. Feigin and E.~Feigin, 
\newblock $q$-characters of the tensor products in 
$\slt$-case, 
math.QA/0201111 (2002).

\bibitem[FFu]{FFu}
B.~L. Feigin and D.~B. Fuchs, 
\newblock Cohomology of some nilpotent subalgebras of the Virasoro and
Kac-Moody algebras, 
\newblock {\em Geometry and Physics},
Essays in honor of I.M.Gelfand on the occasion of his 75th birthday, 
eds. S.Gindikin and I.M.Singer, 209--235, 
North Holland, (1991)


\bibitem[FKLMM1]{FKLMM1}
B.~Feigin, R.~Kedem, S.~Loktev, T.~Miwa and E.~Mukhin,
\newblock Combinatorics of $\slth$ spaces of coinvariants, 
\newblock {\em Transformation Groups} {\bf 6}  (2001) No.1 25--52.

\bibitem[FKLMM2]{FKLMM2}
B.~Feigin, R.~Kedem, S.~Loktev, T.~Miwa and E.~Mukhin,
\newblock Combinatorics of the $\slth$ spaces of Coinvariants:
Loop Heisenberg modules and recursion, 
\newblock math.QA/0009198

\bibitem[FKLMM3]{FKLMM3}
B.~Feigin, R.~Kedem, S.~Loktev, T.~Miwa and E.~Mukhin,
\newblock Combinatorics of the $\slth$ spaces of coinvariants:
Dual functional realization and recursion, 
\newblock math.QA/0012190

\bibitem[FL]{FL}
B.~L. Feigin and S.~Loktev, 
\newblock On generalized Kostka polynomials and 
quantum Verlinde rule,
math.QA/9812093, 
\newblock {\em Amer.~Math.~Sci.~Transl.}
{\bf 194} (1999) 61--79.

\bibitem[FM]{FM}
W.~Fulton and P.~MacPherson. 
\newblock A compactfication of configuration spaces, 
\newblock {\em Ann. Math.} {\bf 139} (1994) 183--225. 

\bibitem[Fo]{Fo}
H.~O.~Foulkes, 
\newblock A survey of some combinatorial aspects of symmetric functions, 
\newblock in Permutations, Gauthier-Villars (1974), Paris.

\bibitem[FS]{FS}
B.~L. Feigin and A.~V. Stoyanovsky,
\newblock Functional models for representations of current algebras and 
  semi-imfinite {Schubert} cells, 
\newblock {\em Funct.~Anal.~and~Its~Appl.}  {\bf 28} (1993) 55--72.

\bibitem[HKKOTY]{HKKOTY}
G.~Hatayama, A.~N.~Kirillov, A.~Kuniba, 
M.~Okado, T.~Takagi and Y.~Yamada, 
\newblock Character formulae of $\widehat{sl}_n$-modules
and inhomogeneous paths, 
math.QA/9802085, 
\newblock {\em Nucl. Phys. B} {\bf 536} (1999) 575--616.

\bibitem[KKMM]{KKMM}    
R.~Kedem, T.~R.~Klassen, B.~M.~McCoy and E.~Melzer, 
\newblock Fermionic sum representations for conformal field theory
    characters, 
\newblock {\em Phys. Lett. B} {\bf 307} (1993) 68--76.

\bibitem[KM]{KM} R.~Kedem and B.~M.~McCoy, 
\newblock Construction of modular branching functions from Bethe's
    equations in the 3-state potts chain,  
\newblock {\em J.  Stat. Phys.} {\bf 71}  (1993) 
875--901.


\bibitem[KR]{KR}
A.~N.~Kirillov and N.~Yu.~Reshetikhin,
\newblock The Bethe Ansatz and the combinatorics of Young tableaux,
\newblock J. Soviet Math. {\bf 41} (1988) 925-955.

\bibitem[KS]{KS}
A.~N.~Kirillov and M.~Shimozono,
\newblock A generalization of the Kostka-Foulkes polynomials,
\newblock  math.QA/9804039,
\newblock J. Algebraic Combin. 15 (2002) 27--69.

\bibitem[KSS]{KSS}
A.~N.~Kirillov, M.~ Shimozono and A.~Schilling, 
\newblock A bijection between Littlewood-Richardson tableaux and rigged
configurations. 
\newblock (Preprint math.CO/9901037).

\bibitem[LS]{LS}
A.~Lascoux and M.P.~Sch\"{u}tzenberger,
\newblock Sur une conjecture de H. O. Foulkes,
\newblock {\em C. R. Acad. Sci. Paris Ser. A-B} {\bf 286} (1978) A323--A324.

\bibitem[Mac]{Mac} I.~G~.Macdonald, 
{\it Symmetric functions and Hall polynomials, 
2nd ed.}, 
Oxford University Press, New York, 1995.

\bibitem[NY]{NY}
A.~Nakayashiki and Y.~Yamada,
\newblock Kostka polynomials and energy functions in solvable lattice models,
\newblock {\em Selecta Math.} {\bf 3} (1997) 547--599.

\bibitem[SS]{SS}
A.~Schilling and M.~Shimozono, 
\newblock Fermionic formulas for level-restricted generalized
Kostka polynomials and coset branching functions, 
math.QA/0001114, 
\newblock {\em Commun. Math. Phys.} {\bf 220} (2001) 105--164, 

\bibitem[SW1]{SW1}
A.~Schilling and  O.~Warnaar, 
\newblock Supernomial coefficients, polynomial identities and $q$-series, 
\newblock Ramanujan J. 2 (1998) 459--494.  

\bibitem[SW2]{SW2}
A.~Schilling and S.~O.~Warnaar,
\newblock Inhomogeneous lattice paths, generalized Kostka polynomials and
$A_{n-1}$ supernomials,
math.QA/9802111,  
\newblock {\em Commun. Math. Phys.} {\bf 202} (1999) 359--401.

\bibitem[TUY]{TUY}
A. Tsuchiya, K. Ueno and Y. Yamada,
\newblock Conformal field theory on the universal family
of stable curves with gauge symmetries,
\newblock {\em Adv. Stud. Pure Math.} {\bf 19} (1989) 459-566.






\end{thebibliography}
\end{document}